\documentclass[11pt]{article}
\usepackage[letterpaper, margin=1in]{geometry}
\usepackage{setspace}
\usepackage{natbib}
 \bibpunct[, ]{(}{)}{,}{a}{}{,}%
 \def\BIBand{and}%
\usepackage{multirow} % Make sure to include this in the preamble

\usepackage{amsmath, amssymb, amsthm}

\newtheorem{theorem}{Theorem}

\newtheorem{corollary}{Corollary}

\numberwithin{equation}{section}  % Use section-based numbering (e.g., 1.1, 1.2)

\usepackage{authblk} % For author affiliations
\usepackage{mathtools, algorithm, algpseudocode, amssymb}
\usepackage{soul}
\usepackage{thm-restate}
\usepackage{enumitem}
\usepackage{tikz}
\usetikzlibrary{trees,shapes.geometric,arrows.meta}
\usepackage{bm}
\usepackage{subcaption}
\usepackage{tabularray}
\usepackage{array}
\usepackage{graphicx} % optional, if needed for superscripts formatting
\usepackage{arydshln} % Required for dashed lines in tables
\usepackage{makecell} % For multi-line cells

\newcommand{\comment}[1]{}
\begin{document}
%%%%%%%%%%%%%%%%

% Enter the full title:
\title{Accelerating Point-Based Value Iteration via Active Sampling of Belief Points and Gaussian Process Regression}

\author[1]{Siqiong Zhou}
\author[1]{Ashif S. Iquebal}
\author[2]{Esma S. Gel}

\affil[1]{School of Computing and Augmented Intelligence, Arizona State University, USA}
\affil[2]{College of Business, University of Nebraska-Lincoln, USA}

\date{} % Removes the default date
\maketitle
\begin{abstract}{%
Partially Observable Markov Decision Processes (POMDPs) are fundamental to decision-making under
uncertainty. We introduce a novel scalable approach to accelerate upper bound estimation in Point-Based Value Iteration (PBVI) algorithms, the leading method to solve large-scale POMDPs. PBVI approximates the value function using a set of belief points rather than the entire continuous belief space and relies on lower and upper bounds for convergence. While lower bounds are straightforward to compute, PVBI requires repeated sawtooth projection operations to approximate the upper bound convex hull, significantly increasing the computational burden although many of these sawtooth projections become redundant as the belief set expands. To address this, we infer the upper bound using the upper confidence bound of a Gaussian Process Regression (GP-UCB) fitted over a subset of the most informative reachable belief points--the ones that exhibit linear independence in some high-dimensional Hilbert space. This approach reduces the number of sawtooth projections by 84.3\% on average without compromising the solution quality. We further establish the theoretical consistency of the proposed GP-UCB estimate of the upper bound and show convergence to the true upper bound convex hull. We implement GP-UCB and test its performance using five benchmark finite-horizon POMDPs, demonstrating its effectiveness in estimating upper bounds and improving PBVI performance. GP-UCB reduces computation time by 30\% to 60\% on smaller problems and up to 99.7\% on larger ones, while achieving the same gaps as the pure sawtooth projection method. 
}
\end{abstract}

% Sample 
%\KEYWORDS{deterministic inventory theory; infinite linear programming duality; 
%  existence of optimal policies; semi-Markov decision process; cyclic schedule}

% Fill in data. If unknown, outcomment the field

%%%%%%%%%%%%%%%%%%%%%%%%%%%%%%%%%%%%%%%%%%%%%%%%%%%%%%%
\section{Introduction}
\label{sec:intro}

A Partially Observable Markov Decision Process (POMDP) extends the classic Markov Decision Process (MDP) by addressing situations where the state of the system is only \textit{ partially observable}~\citep{kaelbling1998planning}. In real-world applications, decision makers often lack perfect information about the current state of the system, which complicates both control and information gathering tasks. This partial observability introduces an inherent uncertainty that must be effectively managed for optimal decision making. The POMDP framework is uniquely suited to this challenge, as it allows decisions to be made based on probabilistic estimates of the state of the system, known as \textit{belief states}.

The significance of solving POMDPs lies in their applicability to a wide range of complex real-world problems. POMDPs have been applied in fields such as reliability and maintenance~\citep{cassandra1998survey, papakonstantinou2014planning, kim2018pomdp, song2022value}, aircraft collision avoidance~\citep{temizer2010collision, bai2012unmanned, mueller2016multi}, and medical decision making~\citep{vozikis2009medical, ayer2012or, zois2016sequential, zhang2022diagnostic}. In each of these domains, the ability to account for uncertain observations and dynamically update belief states based on new information is critical for success.

Although POMDPs offer a powerful theoretical framework, solving them exactly for large-scale problems remains computationally infeasible. This limitation arises from the \textit{curse of dimensionality}, as the size of the reachable belief space grows exponentially with the number of states and observations~\citep{papadimitriou1987complexity, madani1999undecidability}. As a result, approximate methods have been developed to provide practical solutions. Among these, one of the most successful is the Point-Based Value Iteration (PBVI) algorithm~\citep{pineau2003pbvi}, which approximates the value function by focusing on a representative set of belief points rather than computing it over the entire belief space, which is a probability simplex, where each belief $\bm{b}$ is a convex combination of the entire unit vectors corresponding to each state. % (which is continuous). 

The purpose of this paper is to advance the current body of knowledge on solving finite-horizon POMDPs efficiently. In particular, we propose a novel method that accelerates the calculation of upper bounds in PBVI by using a Gaussian Process Upper Confidence Bound (GP-UCB). This method offers a data-efficient approximation strategy to improve the computational tractability of solving POMDPs in complex, large-scale environments. Our contribution is particularly significant in finite-horizon problems, where the dynamically changing value function requires repeated updates to both lower and upper bounds, the latter being computationally prohibitive to estimate. To this end, we employ Gaussian Process Regression (GPR) to learn a model of the upper bound estimates for a subset of the most informative belief points. The trained GPR model is then used to predict the entire upper bound convex hull, reducing the computational complexity while maintaining high-quality approximations. GP-UCB is proposed as a conservative estimate of the upper bound convex hull. The main contributions of our work are as follows. 
\begin{enumerate}
    \item The novel use of the GPR for the upper bound approximation, leveraging the non-parametric approximation capabilities of Gaussian Processes (GP)%to infer the convex hull of upper bounds by fitting the upper bound approximations
    (Section~\ref{sec:ProbDesc}). 
    %for a set of belief points. 
    %GP inference reduces the computational complexity by avoiding redundant backups during the traditional upper bound updates. 
    \item The novel use of active learning for the iterative selection of the most informative belief points, referred to as support beliefs, to train the GPR. % rather than training on the entire belief set.%Second, is the iterative selection of most informative belief points, referred to as support beliefs, using an active learning approach to train the GPR as opposed to training on the entire belief set. 
    Training the GPR on this subset of the belief set reduces the computational cost from $\mathcal{O}(N^3)$ to $\mathcal{O}(Nd^2)$, where $N$ is the number of belief points and $d\ll N$ is the number of support beliefs. 
    \item Theoretical guarantees, backed by formal proofs, showing that the upper bounds generated by our method converge to the true upper bound convex hull, leading to increasingly accurate approximations as the belief set expands.
    %Third, we provide theoretical guarantees for the convergence of the upper bounds generated by our method to the true upper bound convex hull, ensuring that the approximations become increasingly accurate as the belief set expands.
\end{enumerate}

Through extensive numerical experiments, we demonstrate that our method outperforms existing PBVI algorithms in both convergence speed and scalability across a range of complex problem domains and horizon lengths. Our experiments are conducted on well-established test problems, including instances with up to 90 states, 29 actions, and 3 observations. GP-UCB significantly improves computational efficiency in these settings, reducing computation time by 30\%-60\% on smaller problems and up to 99.7\% on larger ones while maintaining the same gaps as the pure sawtooth projection method. 

While POMDP formulations can involve even larger state and action spaces, our experiments exceed the scale of many Operations Research applications (e.g., in medical decision-making and maintenance planning), including recent ones, which typically consider formulations with at most 6 states, 3 actions, and 4~observations~\citep{lin2004hybrid,ayer2012or,liu2022machine,hajjar2023personalized,gong2023partially,li2023optimizing,deep2023partially}. By demonstrating strong performance on significantly larger problems than those addressed in these domains, our approach contributes to the OR literature by providing a scalable and computationally efficient tool for solving practical finite-horizon POMDPs. 

Section~\ref{sec:ProbDesc} provides a description of the main problem that we consider, followed by a detailed explanation of the use of GPR to infer upper bounds in Section~\ref{sec:GPRmethod}. Some readers may find the background given in Appendix~\ref{sec:background} useful before proceeding with these sections since it provides an overview of POMDP solution approaches to date, with a focus on key components of state-of-the-art PBVI algorithms. We demonstrate our effectiveness claims in Section~\ref{sec:perfresults} through the use of an extensive set of numerical experiments on finite-horizon problems described in Section~\ref{sec:NumerExp}. The paper concludes with a summary and comments for future research directions in Section~\ref{sec:Conclusion}.

\section{Problem Description}
\label{sec:ProbDesc}

A POMDP is expressed by a tuple $({\cal S}, {\cal A}, {\cal O}, \Theta, \Omega, R, \bm{b_0})$, where an agent occupies one of the possible states of the system $s \in \cal S$, which cannot be observed directly. The system changes from one state to the next after the agent takes an action $a \in \cal A$. Function~${\Theta: {\cal S}\times {\cal A}\times {\cal S} \rightarrow \left[0,1\right]}$ represents the stochastic state transitions. Specifically, ${\Theta(s,a,s')=P(s'|s,a)}$ denotes the probability transition function of state changes from $s \in \cal S$ to $s' \in \cal S$ after performing the corresponding action $a \in \cal A$. Since the states are not observable directly, the agent makes certain observations $o \in \cal O$, which are imperfect projections of the states. The observation probability is defined by the function~$\Omega: {\cal A} \times {\cal S} \times {\cal O}  \rightarrow \left[0,1\right]$. Specifically, $\Omega(a,s',o)=P(o|a,s')$ is the probability of observing $o \in \cal O$ after taking action $a \in \cal A$, under the true state $s' \in \cal{S}$. In addition, the agent's action results in rewards. Function~$R \in {\cal S}\times{\cal A}\rightarrow \mathbb{R}$ denotes the reward function. Reward~$R(s,a)$ is received after taking action $a \in \cal A$ in state $s \in \cal S$. The planning horizon $T$ specifies the finite time steps during which the agent seeks to maximize its cumulative reward. Let vector ${\bm{b}=(b(1), b(2), \dots, b(|\cal S|))}$ denote the belief state where $b(s)$ is the probability that the true system state is $s \in {\cal S}$. Starting from belief $\bm{b}$ at time~$t$, the agent updates belief to $\bm{b'}$ at time~$t+1$ after executing action $a$ and observing $o$ as
\begin{eqnarray}
b'(s') = \frac{P(o|a,s')}{P(o|b,a)}\sum_{s\in {\cal S}}P(s'|s,a)b(s) = \frac{P(o|a,s')\sum_{s\in {\cal S}}P(s'|s,a)b(s)}{\sum_{s'\in {\cal S}}P(o|a,s')\sum_{s\in {\cal S}}P(s'|s,a)b(s)}~,\forall s \in \cal S~. \label{eq:beliefupdate}
\end{eqnarray}
Naturally, $\sum_{s \in \cal S} b(s) =1$. The initial belief, $\bm{b_0}$, provides the probability distribution of the state at the beginning of the planning horizon. The beliefs derived from an initial belief, via a feasible sequence of actions and observations are called {\it reachable belief points}. 

\begin{figure}[htp]
    \centering
    \includegraphics[width=0.7\textwidth]{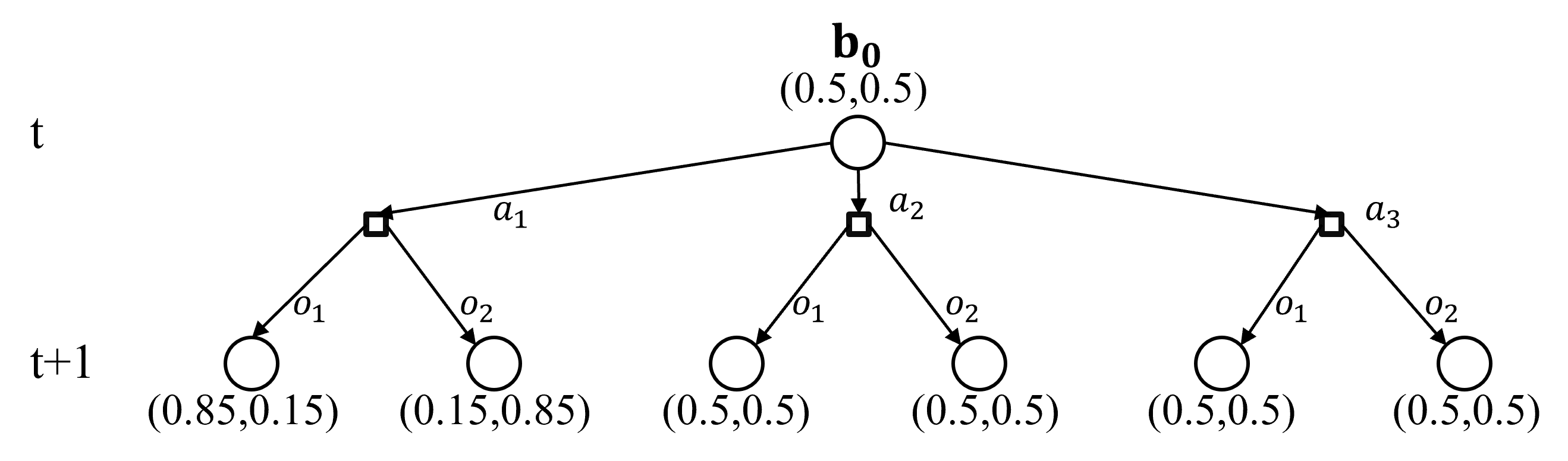}
    \caption{Tree structure of a two-stage tiger problem~\citep{kaelbling1998planning}. Circles represent belief states observed right before taking an action in each stage.}
    \label{fig:POMDPtree}
\end{figure}

Figure~\ref{fig:POMDPtree} shows the belief states that are {\it reachable} from the initial belief state, ${\bm{b_0}=(0.5,0.5)}$ in one stage for the well-known {\it tiger problem} presented by~\citet{kaelbling1998planning}, considering three actions, Listen ($a_1$), Open Left Door ($a_2$), Open Right Door ($a_3$), and two observations resulting from each action, under the assumption that opening a (left or right) door restarts the problem and resets the belief to $(0.5,0.5)$. Appendix~\ref{TigerProblem} contains a detailed description of the tiger problem.

A policy $\pi$ provides the sequence of actions to be taken over the planning horizon as a function of the belief state. The value function represents the expected cumulative reward obtained under the optimal policy, $\pi^*$. For a finite horizon problem, the following backward recursion computes the value function, $V_t(\bm{b})$, as,
\begin{eqnarray}
V_t (\bm{b}) = \underset{a\in {\cal A}}{\max} \left[ \sum_{s\in {\cal S}} b(s)R(s,a)+\sum_{o \in {\cal O}} P(o|\bm{b},a)V_{t+1}^{\pi^*} (\bm{b'}) \right], \forall~ \bm{b} \in \cal{B}_t,~t \in \{0,1,\dots, T-1\}~, \label{valuefunction}
\end{eqnarray}
where $\cal{B}_t$ denotes the $|\cal S|$-dimensional belief space. The optimal policy for belief state $\bm{b}$ at time period $t$ is defined as
\begin{eqnarray}
\pi _t^*(\bm{b}) = \underset{a\in {\cal A}}{\arg \max} \left[ \sum_{s\in {\cal S}} b(s)R(s,a)+\sum_{o \in {\cal O}} P(o|\bm{b},a)V_{t+1}(\bm{b'}) \right]~.
\end{eqnarray}

\citet{sondik1971optimal} showed that the value function is represented exactly by a piecewise-linear and convex function such that the value function for a specific period is represented by a set of $|\cal S|$-dimensional $\alpha$-vectors, $\Gamma_t = \{\alpha^0,\alpha^1,...\}$. That is, 
\begin{eqnarray}
V_t(\bm{b}) = \underset{\alpha \in \Gamma_t}{\max} \sum_{s \in {\cal S}} b(s) \alpha(s)~,~t \in \{0,1,\dots, T-1\}~.
\end{eqnarray}Hence, the value function given in Eqn.~\eqref{valuefunction} can be rewritten as
\begin{eqnarray}
V_t(\bm{b}) = \underset{a\in {\cal A}}{\max} \left[ \sum_{s\in {\cal S}} b(s)R(s,a)+\sum_{o \in {\cal O}} P(o|\bm{b},a) \underset{\alpha_{t+1} \in \Gamma_{t+1}}{\max} \left( \sum_{s'\in {\cal S}} b'(s')\alpha_{t+1}(s') \right) \right]~.\label{eq:lower}
\end{eqnarray}
This representation allows the value function to be computed recursively using a set of $\alpha$-vectors. At each step, the {\it backup operation} updates the set of $\alpha$-vectors, $\Gamma_{t}$ based on the current belief state $\bm{b}$. Specifically, $\Gamma_{t}$ is updated with
\begin{eqnarray}
\Gamma_{t} = \bigcup_{\bm{b} \in \mathcal{B}_t} \text{Backup}(\bm{b}, t)~,\label{eq:gamma}
\end{eqnarray}
where the backup operation for each belief $\bm{b}$ is given by
\begin{eqnarray}
\text{Backup}(\bm{b}, t) =  \underset{g_a^b, \forall a\in {\cal A}}{\arg\max}\left[{\bm{b} \cdot g_a^b}\right]~, \label{eq:backup}
\end{eqnarray}
where $g_a^b$ for action $a$ is computed as
\begin{eqnarray}
g_a^b(s) =
\begin{cases}
R(s,a) + \sum_{o \in \mathcal{O}}  \underset{\alpha_{t+1} \in \Gamma_{t+1}}{\arg\max} \sum_{s \in \mathcal{S}}b(s) \sum_{s' \in \mathcal{S}} P(s'|s,a) P(o | a, s') \alpha_{t+1}(s'), & t < T~, \\
R(s,a), & t = T~.
\end{cases}
\end{eqnarray}

%Recall from the above discussion that 
\paragraph{PBVI Algorithm:} A key idea of PBVI algorithms is sampling a representative set of belief points to approximate the value function. For these sampled belief points, both a lower bound $\underline{V_t}$ and an upper bound $\overline{V_t}$ of the value function are maintained. The gap between the lower and upper bounds at reachable belief points affects the value function approximation. Thus, the PBVI algorithm implemented by \citet{walraven2019point} selects the belief points with the largest gap, focusing on areas that can most effectively improve the value function approximation. The approach prioritizes the sampling of belief points where the gap is maximum in the time step~$t+1$. 

Algorithm~\ref{algo:PBVI} outlines the PBVI algorithm incorporating lower and upper bound calculations. Following \citet{walraven2019point}, the initial belief set includes corner beliefs, which are the belief points where all probability mass is concentrated on a single state. Mathematically, these are the unit vectors $\bm{w}_s$, where $\bm{w}_s$ is a vector with $1$ in the $s^{th}$ entry and $0$ elsewhere. These points define the boundaries of the belief space, and hence are useful for initializing the PBVI and defining the reachable belief set.
% \begin{equation}
% b(s) =
% \begin{cases}
%     1, & \text{if } s = s^* \text{ for some } s^* \in \mathcal{S}, \\
%     0, & \text{otherwise}.
% \end{cases}
% \label{eq:corner_belief}
% \end{equation}

\begin{algorithm}[h]
\small
\caption{Point-based Value Iteration (PBVI) Algorithm, generalized from \citet{walraven2019point, spaan2005perseus,lovejoy1991computationally}}
\label{algo:PBVI}
\begin{algorithmic}[1] % The [1] option enables line numbering
\State \textbf{Input:} {POMDP model, initial belief point $\bm{b_0}$, initial belief set $\cal{B}_t$ containing all corner beliefs for $t \in \{0,1, \dots, T-1\}$, and convergence threshold, $\epsilon$}
\State \textbf{Output:} {Lower bound $\underline{V_t}(\bm{b})$ and upper bound $\overline{V_t}(\bm{b})$ of value function for $\forall \bm{b} \in \cal{B}_{t}$ for $t \in \{0,1, \dots, T-1\}$}
\State Initialize $\underline{V_t}(\bm{b})$ and $\overline{V_t}(\bm{b})$ for $\forall \bm{b} \in \cal{B}_{t}$ for $t \in \{0,1, \dots, T-1\}$ using Eqn.~\eqref{eq:lower} and Eqn.~\eqref{eq:upperbound}
\While{{$\overline{V_0}(\bm{b_0}) - \underline{V_0}(\bm{b_0}) > \epsilon$} }
    \State Choose a sampling method: Max-Gap, Random or Fixed-Grid
    \If{Sampling method = ``Max-Gap"} 
        \State {Set $\bm{b} = \bm{b_0}$}
        \For{{$t = 0$ to $T-2$}}
            \State Choose the action $a \gets \arg\max_{a \in {\cal A}} \left(  R(a) \cdot \bm{b} +  \sum_{o \in {\cal O}} P(o|\bm{b},a)\cdot \overline{V_{t+1}} (\bm{b'}) \right)$
            \State Choose the observation $o \gets \arg\max_{o \in {\cal O}} \left(  \overline{V_{t+1}}(\bm{b}) - \underline{V_{t+1}}(\bm{b}) \right)$
            \State Sample new belief, $\bm{b_{\text{new}}}$ using the selected $a$, $o$ with Eqn.~\eqref{eq:beliefupdate} and add to $\cal{B}_{t+1}$
            \State {Set $\bm{b} = \bm{b_{\text{new}}}$}
        \EndFor \hspace{0.5cm} {\textit{(The max-gap method from \citet{walraven2019point})}}
    \ElsIf{Sampling method = ``Random"}
        \For{{$t = 1$ to $T-1$}}
            \State Sample new belief, $\bm{b_{\text{new}}}$, uniformly from belief space and add to $\cal{B}_{t}$, e.g., random sampling method proposed by \citet{spaan2005perseus}
        \EndFor
    \ElsIf{Sampling method = ``Fixed-Grid"}
        \State Use fixed belief grid for $\cal{B}_{t}$ for $t \in \{0,1, \dots, T-1\}$, e.g., fixed grid proposed by \citet{lovejoy1991computationally}
    \EndIf
    \For{{$t = T-1$ to $0$}}
        \State Update $\underline{V_t}(\bm{b})$ for $\forall \bm{b} \in \cal{B}_{t}$ using Eqn.~\eqref{eq:lower}
        \State Update $\overline{V_t}(\bm{b})$ for $\forall \bm{b} \in \cal{B}_{t}$ with Eqn.~\eqref{eq:upperbound} using an approximation method, e.g., sawtooth
        \State Prune dominated $\alpha$-vectors
    \EndFor
\EndWhile
\end{algorithmic} % The [1] option enables line numbering
\end{algorithm}

We now discuss the computation of $\underline{V_t}(\bm{b})$ and $\overline{V_t}(\bm{b})$ for a belief point, $\bm{b} \in \cal{B}_t$. While the former can be easily obtained by taking the best $\alpha$-vector at any belief point~\citep{hauskrecht2000value}, obtaining the upper bound is computationally more challenging. In finite-horizon POMDPs, the challenge is even greater because the value function evolves dynamically over time \citep{smallwood1973optimal}. Unlike in infinite-horizon settings—where a stationary policy optimizes the value function indefinitely—finite-horizon problems require recomputing the value function at each time step, as the number of remaining decisions directly influences both immediate and future rewards~\citep{pineau2003pbvi}. This makes the calculation of upper bounds particularly expensive, since upper bounds must be recomputed at each time step to reflect the evolving decision horizon. Determining an exact upper bound requires optimistically evaluating all possible future outcomes, significantly increasing computational cost~\citep{walraven2019point, smith2005point}. %This is the main reason that we test our approach on finite-horizon problems in this paper. 

\citet{lovejoy1991computationally} presented some of the earliest implementation of upper bounds using linear interpolation method using a grid of belief points obtained via Freudenthal triangulation. Subsequently, \cite{hauskrecht2000value} presented an upper bound based on the convex hull projection while considering the upper bound improvement with incremental addition of new belief points. Since the value function is piecewise convex, the upper bound for a new belief point is obtained by projecting the belief point to the convex hull obtained by the existing belief-upper bound value pairs. For the remainder of the paper, we refer to this as {\it convex hull}. Linear programming is used to identify the best convex hull by minimizing the linear combination of the existing upper bounds. Given the repetitive nature of updating upper bounds, the computational cost of executing a large number of linear programming steps quickly becomes intractable.

\cite{hauskrecht2000value} proposed an interpolation approach for efficient approximation of the convex hull, referred to as {\it sawtooth}, to reduce computational complexity. Given an arbitrary belief $\bm{b} \in \cal{B}$ and corner beliefs, the sawtooth projection $\bar{v}(\cdot)$ provides the upper bound approximation for any new belief point. The detailed description of the sawtooth projection algorithm can be found in Appendix~\ref{sec:sawtooth}.

\begin{figure}[htp]
    \centering
    \includegraphics[width=1\textwidth]{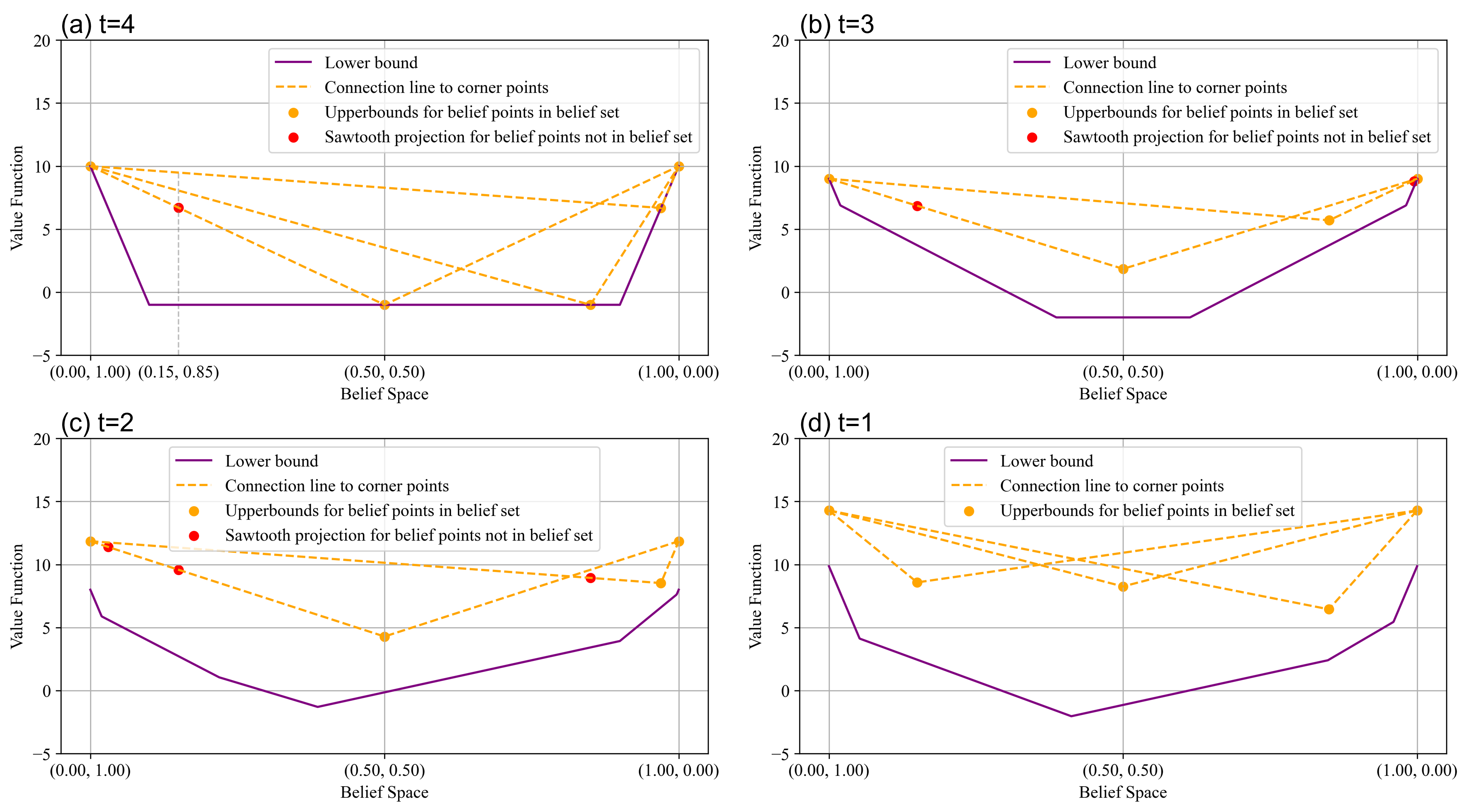}
    \caption{Approximation of the upper bound using sawtooth projection across different time stages at iteration~2 in the tiger problem.}
    \label{fig:SawtoothTiger}
\end{figure}

To further demonstrate the sawtooth projection, consider the tiger problem again. In each subfigure of Figure~\ref{fig:SawtoothTiger}, the orange dots represent belief points in the current belief set $\cal{B}$ with known upper bounds. In the tiger problem, corner beliefs correspond to states where the tiger is surely behind either the left or the right door (i.e., belief states $(1.00, 0.00)$ or $(0.00, 1.00)$). Connecting non-corner orange belief points with these corner beliefs forms downward-pointing triangles, illustrated by the orange dashed lines. When a new belief $\bm{b'}$ outside $\cal{B}$ is encountered, its $\bar{v}(\bm{b'})$ is estimated by projecting it onto these connecting lines, selecting the projection with the smallest value (marked with red dots). For example, when updating the upper bound for belief points at $t=3$, such as $(0.5,0.5)$, Eqn.~\eqref{eq:upperbound} is applied, requiring $\bar{v}(\cdot)$ for reachable beliefs at the subsequent stage, $t=4$. In this case, the reachable belief points are $(0.5,0.5)$, $(0.15,0.85)$, and $(0.85,0.15)$. Figure~\ref{fig:SawtoothTiger}(a) shows that $(0.15,0.85)$ is not part of the current belief set, so $\bar{v}[(0.15,0.85)]$ is approximated using a sawtooth projection. By projecting $(0.15,0.85)$ onto the dashed orange lines, the smallest projection value is selected, shown as the red point in Figure~\ref{fig:SawtoothTiger}(a). For other reachable belief points whose upper bounds are known (orange points), the values are used as $\bar{v}(\cdot)$ directly. This process repeats as we move back to $t=2$, $t=1$, and $t=0$, applying sawtooth projections to any reachable points with unknown upper bounds, as depicted by the red points in Figures~\ref{fig:SawtoothTiger}(b)-(d).

The sawtooth projection method offers a trade-off between computational efficiency and accuracy. Using geometric properties, it efficiently approximates the convex hull, avoiding the extensive recalculations required by heuristic-based approaches such as HSVI~\citep{smith2005point} and SARSOP~\citep{kurniawati2008sarsop}. However, the sawtooth projection still presents significant computational complexity due to the repeated updates required to refine the upper bound each time a new belief point is added to the belief set. Initially, expanding the belief set results in significant improvements in the upper bound. But as the belief set grows, the improvement diminishes, although the number of updates, and consequently the computational complexity, continues to increase~\citep{hauskrecht2000value}.

\section{Upper Bound Prediction using Gaussian Process Regression}
\label{sec:GPRmethod}

Here, we introduce an alternative approach that approximates the upper bound convex hull by fitting a GPR to the sawtooth projections of only a subset of belief points that we refer to as the {\it support beliefs}. As noted in \cite{hauskrecht2000value}, the change in upper bounds as the belief set expands is only marginal. As a result, an  ``informative'' subset of the belief set is typically sufficient to model the convex hull. %An outline of the proposed approach is as follows. 

Recall that under the PBVI algorithm presented by \citet{walraven2019point}, at any stage of upper bound updates for a set of belief points, we execute the backup operation as shown in Eqn.~\eqref{eq:upperbound}, which involves computing the sawtooth projection for all the reachable beliefs. Instead, we compute the sawtooth projection only for the support beliefs and predict for the rest of the belief points using a GPR. The support belief set is initialized with $|\mathcal{S}|+1$ random belief points and iteratively expanded using an approximate linear dependence criterion discussed below. Finally, to ensure that our approximation provides a conservative estimate of the convex hull, we consider the upper confidence bound of the GPR during the inference stage.

GPR learns an interpolation of the convex hull using a set of $m$ arbitrary reachable belief points and the corresponding sawtooth projections $\{\bm{b}_i,\overline{v}\}_{i=1}^m$. In this framework, the sawtooth projections $\overline{v}(\bm{b}_i)$ represent the noisy estimates of the convex hull, $h(\bm{b}_i)$, such that $\overline{v}(\bm{b}_i) = h(\bm{b}_i)+\epsilon_i$, where ${\epsilon_i \sim \mathcal{N}(\mathbf{0},\sigma^2_\epsilon)}$ is assumed to be normally distributed noise with variance $\sigma^2_\epsilon$, representing the variability in the overestimation of the sawtooth projections. Hence, the sawtooth projection of the training set follows a multivariate Gaussian distribution, i.e., $\overline{\bm{v}}_m \equiv [\overline{v}(\bm{b}_1),\ldots, \overline{v}(\bm{b}_m)]^T \sim \mathcal{N} \left( \mathbf{0}, \mathbf{K}(\boldsymbol{B}_m,\boldsymbol{B}_m) + \sigma^2_\epsilon \mathbf{I} \right)$, where $\boldsymbol{B}_m = {\{\bm{b}_i\}_{i=1}^m}$ is the support belief set and $\mathbf{K}(\boldsymbol{B}_m,\boldsymbol{B}_m)$ denotes the covariance matrix. For any new belief point $\bm{b}_{k}$, the joint distribution of known upper bounds $\overline{\bm{v}}_m$ and estimated upper bound $\widehat{v}(\bm{b}_{k})$ is Gaussian,~i.e. 
\begin{eqnarray}
\begin{bmatrix}
\overline{\bm{v}}_m\\
\widehat{v}(\bm{b}_{k})
\end{bmatrix} \sim \mathcal{N} \left( \mathbf{0}, \begin{bmatrix}
\mathbf{K}(\boldsymbol{B}_m,\boldsymbol{B}_m)+\sigma^2_\epsilon \mathbf{I} & \mathbf{K}(\boldsymbol{B}_m,\bm{b}_{k})\\
\mathbf{K}(\bm{b}_{k},\boldsymbol{B}_m) & \mathbf{K}(\bm{b}_{k},\bm{b}_{k})
\end{bmatrix} \right)
\end{eqnarray}
where $\overline{\bm{v}}_m$ represents the upper bound at the current belief set $\boldsymbol{B}_m$, and $\widehat{v}(\bm{b}_{k})$ is the estimated upper bound at the new belief point $\bm{b}_{k}$.Conditioning on the joint distribution yields the predicted mean value $\mu(\bm{b}_{k})$ of the convex hull at $\bm{b}_{k}$ with the corresponding variance $\sigma^2(\bm{b}_{k})$ at a new belief point $\bm{b}_{k}$ as
\begin{align}
 \mu(\bm{b}_{k}) &= \mathbf{K}(\boldsymbol{B}_m,\bm{b}_{k})^T \left[ \mathbf{K}(\boldsymbol{B}_m,\boldsymbol{B}_m) + \sigma^2_\epsilon \mathbf{I} \right]^{-1} \overline{\bm{v}}_m~, \label{eq:GPRMean} \\
\sigma^2(\bm{b}_{k}) &= \mathbf{K}(\bm{b}_{k},\bm{b}_{k}) - \mathbf{K}(\boldsymbol{B}_m,\bm{b}_{k})^T \left[ \mathbf{K}(\boldsymbol{B}_m,\boldsymbol{B}_m)+ \sigma^2_\epsilon \mathbf{I} \right]^{-1} \mathbf{K}(\boldsymbol{B}_m,\bm{b}_{k})~\label{eq:GPRVar}.
\end{align}

% \subsection{Update the Gaussian Model}
% \label{sec:UpdateGP}
Let us define $\bm{\beta} = \left[ \mathbf{K}(\boldsymbol{B}_m,\boldsymbol{B}_m) + \sigma^2_\epsilon \mathbf{I} \right]^{-1} \overline{\bm{v}}_m$ such that Eqn.~\eqref{eq:GPRMean} is rewritten as the following linear combination of kernel functions $k(\bm{b}_i, \bm{b}_{k})  = [\mathbf{K}(\boldsymbol{B}_m, \bm{b}_{k})]_i$, 
\begin{equation}
    \mu(\bm{b}_{k}) = \sum_{i = 1}^{m}\beta_i k(\bm{b}_i, \bm{b}_{k}) = \sum_{i = 1}^{m} \beta_i \langle \phi(\bm{b}_i), \phi(\bm{b}_{k})\rangle~,
\end{equation} where $\phi$ is a mapping of the points in the belief space to Hilbert space. This form of the mean GPR predictor shows that it is, in fact, a linear predictor of the convex hull in the Hilbert space. As such, if an arbitrary belief point $\bm{b}_{k}$ satisfies $\phi(\bm{b}_{k}) = \sum_{i = 1}^{m}\beta_i \phi(\bm{b}_i)$, then the mean GPR prediction corresponding to $\bm{b}_{k}$ could be inferred directly from the beliefs $\bm{b}_1, \ldots, \bm{b}_{m}$. However, unless for the cases when the belief space is low dimensional or the belief $\bm{b}_{k} = \bm{b}_j, j \leq m$, $\phi(\bm{b}_{k})$ will be linearly independent of $\{\phi(\bm{b}_i)\}_{i=1}^{m}$. Although strong linear dependency does not hold, \citet{engel2004kernel} showed that it is possible to have an approximate linear dependency. To test whether a new belief $\bm{b}_{k}$ satisfies the approximately linear dependent or ALD criterion, we define 
\begin{equation}
\delta \; \stackrel{\text{def}}{=} \; \min_{\pmb{a}}{\Bigg\|\sum_{j=1}^{m}a_j\pmb{\phi}(\bm{b}_j)-\pmb{\phi}(\bm{b}_{k}) \Bigg\|}^2 = \; k(\bm{b}_{k}, \bm{b}_{k})-{\mathbf{K}}(\boldsymbol{B}_m, \bm{b}_{k})^\top\pmb{a}\; \leq \; \nu~, \label{eq:ALD}
\end{equation}
\noindent where $\pmb{a}\equiv [a_1, a_2, \ldots, a_m]^T = {\mathbf{K}}^{-1}(\boldsymbol{B}_m, \boldsymbol{B}_m){\mathbf{K}}(\boldsymbol{B}_m, \bm{b}_{k})$, and $\nu$ is a threshold that determines the strength of ALD. Here, $\mathbf{K}^{-1}(\boldsymbol{B}_m, \boldsymbol{B}_m)$ is the inverse of the kernel matrix for the belief set $\boldsymbol{B}_{m}$. If $\delta \leq \nu$, then $\bm{b}_{k}$ is considered ALD on the support belief $\boldsymbol{B}_{m}$. In this case, the support belief set remains unchanged. However, if $\delta > \nu$, $\bm{b}_{k}$ is not linearly dependent on $\boldsymbol{B}_{m}$, in which case the support belief $\boldsymbol{B}_m$ is expanded to include the new belief point $\bm{b}_{k}$, and the corresponding kernel matrices are updated accordingly. 

Note, however, that one needs to repeat the GPR fitting every time a new belief point $\bm{b}_{k}$ is added to the belief set. This is because as new belief points are added,  the sawtooth projections $\overline{\bm{v}}_m$ of existing beliefs in the training set $\boldsymbol{B}_{m}$ get tighter, leading to a more accurate upper bound approximation of the value function. However, these updated projections are not directly available and must be recomputed. 

To reduce computational load, we update the GPR fit using partially revised sawtooth projections, which are then used to obtain the updated mean prediction in Eqn.~\eqref{eq:GPRMean}, rather than retraining from scratch when upper bound changes are small. Specifically, during initial iterations, when upper bounds change significantly, the GPR model is refitted at each iteration to ensure accuracy. However, as the upper bounds stabilize in later iterations, a single belief point is randomly selected from $\boldsymbol{B}_{m}$ in each iteration, and its sawtooth projection is recomputed to update the GPR fit in Eqn.~\eqref{eq:GPRMean}. Full updates of all sawtooth projections $\overline{\bm{v}}_m$ become unnecessary because upper bound improvements tend to stabilize after a few iterations, as also noted by \citet{hauskrecht2000value} and \citet{smith2012point}.

\begin{figure}[!b]
    \centering
    \includegraphics[width=1\textwidth]{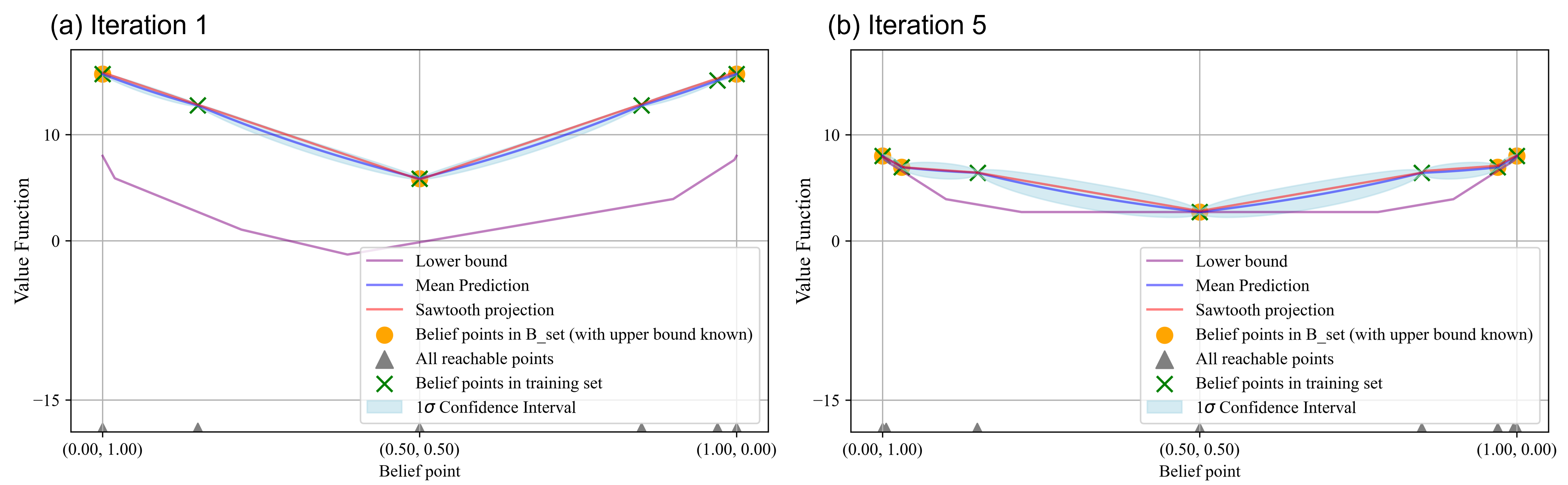}
    \caption{Upper bound estimation using GPR and sawtooth projection. } 
    \label{fig:GPRetimation}
\end{figure}

Figure~\ref{fig:GPRetimation} provides a visualization of estimating upper bounds using GPR in the tiger problem example. Figure~\ref{fig:GPRetimation}(a) shows the GPR model fit in iteration~1, and Figure~\ref{fig:GPRetimation}(b) presents the fitting at iteration~5. At iteration~1, the support belief set includes three belief points with known upper bounds, depicted as orange dots, and three additional points selected based on the ALD criterion, with their respective sawtooth projections shown as green x markers. The GPR model is trained using these six belief points. The resulting convex hull inferred from the GPR model is shown in a solid blue line together with the $1\sigma$-confidence bound. For comparison, the solid red line shows the upper bound approximation obtained from the sawtooth projection. As the belief set expands over subsequent iterations, the GPR model is updated with the new belief points using the ALD criterion. At iteration~5, the GPR model is trained by incorporating two additional reachable belief points alongside five existing beliefs with known upper bounds. The updated fitting of the GPR model reflecting the refined upper bound approximation, is depicted in Figure~\ref{fig:GPRetimation}(b).

\subsection{The Use of Gaussian Process Upper Confidence Bound}
\label{sec:UCB}

\begin{figure}[!htp]
    \centering
    \includegraphics[width=1\textwidth]{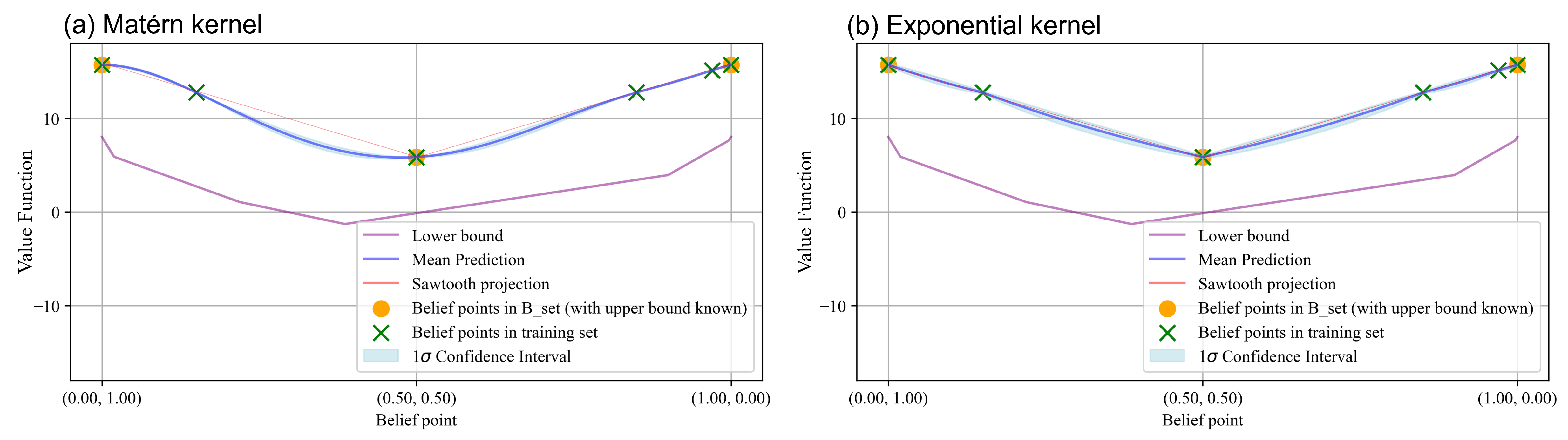}
    \caption{Upper bound estimation using Mat\'ern and exponential kernels in the tiger problem iteration~1. }
    \label{fig:GPRkernel}
\end{figure}

The last issue we address in this work is that of GP mean reversal and its impact on the convex hull prediction. The predicted mean tends to collapse toward the prior mean in regions with no or fewer training points. Since the prior mean is assumed to be 0 in the current case, the predicted convex hull would have a tendency to sag toward the zero line. The mean reversal effect is more pronounced with smoother kernel functions such as squared exponential and Mat\'ern, but less so with exponential. An example comparison is shown in Figure~\ref{fig:GPRkernel}. 

Due to the mean reversal effect, the predicted convex hull may underestimate the true one. Hence, to overcome this limitation, we apply the Upper Confidence Bound (UCB) of the predicted convex hull as the upper bound $\overline{V}_t(\bm{b})$. For any belief point $\bm{b}$, the upper confidence bound $\overline{V}_{UCB}(\bm{b}) = \mu(\bm{b}) + \eta \sigma(\bm{b})$ provides a conservative estimate of the convex hull ${h}(\bm{b})$, where $\eta$ is a confidence parameter. In the following, we theoretically show that the proposed upper bound $\overline{V}_{UCB}(\cdot)$ is a Probably Approximately Correct (PAC) estimate of the convex hull. 

\begin{theorem}
\label{theo:convergeCH}
Denoting the convex hull upper bound for belief point $\bm{b}$ as ${h}(\bm{b})$, we have 
\[
P\left( |{h}(\bm{b})-\mu (\bm{b})| \leq \eta \sigma(\bm{b}) \right) \geq 1-\delta \text{ for some } \delta \in (0,1)~.
\]
In other words, as more belief points are sampled, $\sigma(\bm{b})$ reduces and the mean prediction $\mu (\bm{b})$ approaches the true convex hull $h(\bm{b})$.
\end{theorem}

\begin{proof}{Proof.}
Conditioned on a given set of belief points \( \boldsymbol{B} = {\{\bm{b}_i\}_{i=1}^m} \) and their corresponding sawtooth projections for upper bounds $\overline{\bm{v}}(\bm{b}_i)$, $i=1,\dots,m$, the posterior distribution of the estimated upper bound using GPR is $\widehat{v}(\bm{b}) \sim \mathcal{N}(\mu(\bm{b}), \sigma(\bm{b}))$,
where $\mu(\bm{b})$ is the mean prediction and $ \sigma(\bm{b})$is the standard deviation.

As discussed in \cite{hauskrecht2000value}, the sawtooth projection \( \bar{v}(\bm{b}) \) overestimates the convex hull upper bound \( h(\bm{b}) \), ensuring \( \bar{v}(\bm{b}) \geq h(\bm{b}) \). Since the GPR model is trained on sawtooth projections, the posterior mean $\mu(\bm{b})$ progressively approximates $h(\bm{b})$ as additional belief points are sampled. Using the properties of Gaussian distribution, the normalized deviation of \( h(\bm{b}) \) from \( \mu(\bm{b}) \), relative to the standard deviation \( \sigma(\bm{b}) \), satisfies
$$P\left( \frac{|h(\bm{b}) - \mu(\bm{b})|}{\sigma(\bm{b})}\leq \eta \right) \geq \frac{1}{\sqrt{2 \pi}} e^{-\eta^2/2}~.$$
Rearranging this inequality, we obtain
$P\left( |h(\bm{b}) - \mu(\bm{b})| \leq \eta \sigma(\bm{b}) \right) \geq 1 - \delta~,$
where ${\delta = 1 - \frac{1}{\sqrt{2 \pi}} e^{-\eta^2/2}}$.
\qed
\end{proof}

Theorem~\ref{theo:convergeCH} indicates that as more belief points are added, the uncertainty in the convex hull prediction, \( \sigma(\bm{b}) \) decreases, reflecting greater confidence in the GPR prediction. Consequently, the posterior mean \( \mu(\bm{b}) \) converges to the convex hull upper bound \( h(\bm{b}) \).

\begin{corollary}
\label{coro:UCBestimate}
As $\mu (\bm{b})$ approaches $h(\bm{b})$, the proposed upper confidence bound $\overline{V}_{UCB}(\bm{b})$ becomes a conservative estimate of true upper bound, i.e.
\[
P\left( h(\bm{b}) \leq \overline{V}_{UCB} (\bm{b}) \right) \geq  \frac{1}{\sqrt{2 \pi}} e^{-\eta^2/2}~.
\]
\end{corollary}
\begin{proof}{Proof.}
Consider the definition of the proposed upper confidence bound \( \overline{V}_{UCB}(\bm{b}) \):
$\overline{V}_{UCB}(\bm{b}) = \mu(\bm{b}) + \eta \sigma(\bm{b}).$
Theorem~\ref{theo:convergeCH} ensures that \( \mu(\bm{b}) \) is close to \( h(\bm{b}) \) with high probability. Therefore, the probability that the  upper bound convex hull \( h(\bm{b}) \) lies below the proposed upper confidence bound \( \overline{V}_{UCB}(\bm{b}) \) is
$P\left( h(\bm{b}) \leq \overline{V}_{UCB}(\bm{b}) \right) = P\left( h(\bm{b}) - \mu(\bm{b}) \leq \eta \sigma(\bm{b}) \right).$
From Theorem~\ref{theo:convergeCH}, we know that 
$P\left( |h(\bm{b}) - \mu(\bm{b})| \leq \eta \sigma(\bm{b}) \right) \geq 1 - \delta.$ Substituting ${\delta = 1 - \frac{1}{\sqrt{2 \pi}} e^{-\eta^2/2}}$, we get
$$P\left( h(\bm{b}) \leq \overline{V}_{UCB}(\bm{b}) \right) \geq \frac{1}{\sqrt{2\pi}} e^{-\eta^2/2}~.$$
Thus, as \( \mu(\bm{b}) \) converges to \( h(\bm{b}) \) with more belief points sampled and \( \sigma(\bm{b}) \) decreases, the proposed upper confidence bound becomes a probabilistically conservative estimate of the true upper bound \( h(\bm{b}) \). \qed
\end{proof}

\begin{figure}[htp]
    \centering
\includegraphics[width=1\textwidth]{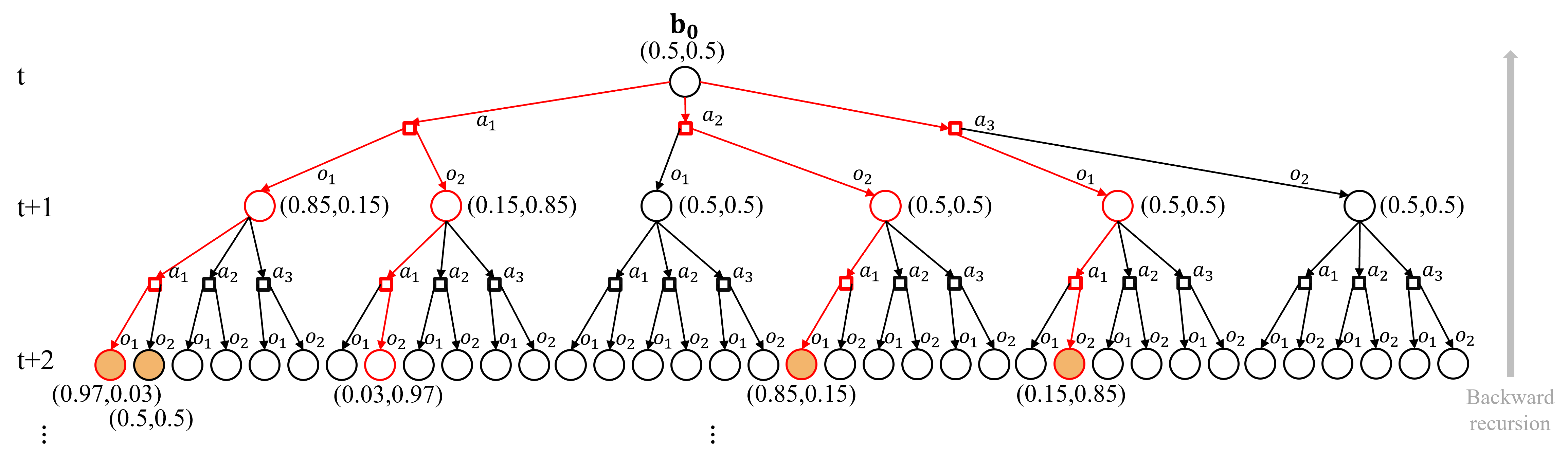}
\caption{Estimation of upper bounds with the GPR model in the tiger problem. }
    \label{fig:TreeStructure-GPR}
\end{figure}

Figure~\ref{fig:TreeStructure-GPR} shows how the GPR model estimates the upper bounds for the tiger problem. Starting from $(0.5,0.5)$, the belief points reachable by $t+2$ are depicted as circles, with sampled beliefs shown on the red routes. The upper bounds of these sampled beliefs are updated using backward recursion. At $t+2$, the GPR model is constructed using the most informative belief points in the belief set (shown with orange-shaded circles). In this example, one of the belief points at $t+2$ (shown with a red edge but not orange-shaded) is not selected to train the GPR due to its lack of information gain. In addition, note that the belief point $(0.5,0.5)$, which is not a belief point originally sampled, is selected to train the GPR model based on the ALD criterion, indicating it is an informative belief point. The trained GPR model is then used to estimate the upper bounds of other reachable belief points at $t+2$ using the GP-UCB. Finally, we estimate upper bounds for belief points at $t+1$ using the backward recursion given by  Eqn.~\eqref{eq:upperbound} in Appendix B. 

Algorithm~\ref{algo:GP_UCB_POMDP} outlines the steps for approximating an upper bound to the value function using our proposed GP-UCB approach.

%\esg{SIQIONG, IN THE ALGORITHM, HOW DO WE INITIALIZE $\Gamma-t$? the term backup operation may not be clear, provide Equation number in the algo. Why did you choose 5 iterations for the if statement? on line 13, why is the training set updated? i cannot see it in the algo? We should say ``updated upper bounds of the points in training set'' Line 15 says randomly select a point b in Bt, but do you use other methods to update the GPR model? On line 20, is that the support beliefs or the training set (which is subset of support beliefs) that you are referring to? makes more sense that you are referring to the set of support beliefs (calligraphic), but then I do not see which line support belief set is updated?-- on line 5, what does Expand belief set $\cal B_t$ mean? How is it done, or is that the header for the for statement? Try to provide equations for update operations, backward recursion, etc. }

\begin{algorithm}[h]
\small
    \caption{GP-UCB Algorithm}
    \label{algo:GP_UCB_POMDP}
    \begin{algorithmic}[1] \setlength{\baselineskip}{0.9\baselineskip} % Reduces line spacing
    \State \textbf{Input:} Planning horizon $T$, initial belief $\bm{b}_0$, threshold $\epsilon$
    \State \textbf{Output:} $\overline{V_0}(\bm{b}_0)$, $\underline{V_0}(\bm{b}_0)$, and $\tt{gap}$ for initial belief $\bm{b}_0$
    
    \State \textbf{Initialization:} For $t \in \{0,1, \dots, T-1\}$, initialize $\cal{B}_t$ to include all corner points and $\bm{b}_0$, initialize training set $\boldsymbol{B}_t$ to include all points in $\cal{B}_t$, and initialize $\Gamma_t$ using  Eqn.~\eqref{eq:gamma}.
    
    \While{running time $<$ time limit and $\tt{gap} > \epsilon$}
        \State Expand $\cal{B}_t$ for $t \in \{0,1, \dots, T-1\}$ using an sampling method, e.g., max-gap (lines 6-13 in Algo.~\ref{algo:PBVI}), random (lines 14-17 in Algo.~\ref{algo:PBVI}), fixed-grid (lines 18-20 in Algo.~\ref{algo:PBVI})
        \For{$t = T-1$ to $0$}
            \State Update $\Gamma_t$ for $\cal{B}_t$ using Eqn.~\eqref{eq:gamma}
            \If{first iteration}
                \State Initialize GPR model for $t$ with beliefs in $\boldsymbol{B}_t$
            \Else
                \If{initial iterations \textbf{or} $\tt{gap}$ change$> 100\cdot\epsilon$ \textbf{or} periodic check iteration}
                    \State Update all upper bounds for beliefs in $\boldsymbol{B}_t$ using sawtooth projection
                    \State Fit the GPR model with updated upper bounds of the beliefs in $\boldsymbol{B}_t$
                \Else
                    \State Randomly select a belief $\bm{b}$ in $\boldsymbol{B}_t$
                    \State Compute the updated sawtooth projection for $\bm{b}$
                    \State Update the GPR fit using revised sawtooth projections
                \EndIf
            \EndIf
            \If{$\cal{B}_t$ expanded with new points = True}
                \State Predict the covariance for new reachable beliefs using the GPR model
                \If{Covariance $\geq$ ALD threshold}
                    \State Expand $\boldsymbol{B}_t$ with the new reachable belief
                    \State Update GPR model using the expanded beliefs in $\boldsymbol{B}_t$
                \EndIf
            \EndIf
            \State Update $\overline{V_t}(\bm{b})$ for $\forall \bm{b} \in \cal{B}_{t}$  using Eqn.~\eqref{eq:upperbound} with GPR model predictions
        \EndFor
        \State Compute $\overline{V_0}(\bm{b}_0)$ and $\underline{V_0}(\bm{b}_0)$ for belief $\bm{b}_0$ and update ${\tt gap} = \overline{V_0}(\bm{b}_0) - \underline{V_0}(\bm{b}_0)$
    \EndWhile
    \State \Return $\overline{V_0}(\bm{b}_0)$, $\underline{V_0}(\bm{b}_0)$, and $\tt{gap}$ for the $\bm{b}_0$
    \end{algorithmic}
\end{algorithm}

% Algorithm~\ref{algo:AUCB} outlines the steps for approximating an upper bound to the value function using our proposed GP-UCB approach.

% \begin{algorithm}[!t]
%     \caption{GP-UCB confidence bound}
%     \label{algo:AUCB}
%     \begin{algorithmic}[1] % Enables line numbering
%     \State \textbf{Input:} POMDP model, a newly sampled belief point $b_{new}$, initialized training set $\{B_{train}, \mathcal{ST}_{train}\}$ with a subset of belief points and corresponding upper bounds
%     \State \textbf{Output:} upper bound $\overline{V}$ of value function
%     \If{Refit criteria = True}
%         \State Update upper bounds $\mathcal{ST}_{train}$ for all belief points in $B_{train}$ by using sawtooth projection.
%         \State Fit the training set $\{B_{train}, \mathcal{ST}_{train}\}$ with Gaussian process regression.
%     \Else
%         \State Randomly pick one belief point $b \in B_{train}$ and update the corresponding $\mathcal{ST}_{train}$ with a sawtooth projection.
%         \State Update the predicting vector of the Gaussian model with the updated training set.
%     \EndIf
%     \If{Add $b_{new}$ to the training set = True}
%         \State Update the Gaussian model with the kernel recursive least-squares algorithm.
%         \State Expand the training set $\{B_{train}, \mathcal{ST}_{train}\}$ with $b_{new}$.
%     \EndIf
%     \State Update $\overline{V}$ with the upper confidence bound from Gaussian prediction.
%     \State \Return{$\overline{V}$}
%     \end{algorithmic}
% \end{algorithm}

\section{Numerical Experiments}
\label{sec:NumerExp}

In this section, we present the experimental evaluation of our approach using five commonly used test problems from pomdp.org~\citep{pomdporg}. These examples are selected to compare the performance of the algorithm on a range of problem sized, including varying numbers of states, actions, observations, and horizon length. 

Table~\ref{ExampleProperty} summarizes the environmental variables for the five test problems used in our experiments. Since we focus on the finite planning horizon, we consider undiscounted versions of the problems (i.e., discount factors are set to $1$). Additionally, we set the initial belief points for all examples as the center point, i.e., $(1/|{\cal S}|, 1/|{\cal S}|, ..., 1/|{\cal S}|)$, which indicates that all states are equally likely at $t=0$. To ensure robustness and account for variability due to the random selection of belief point in updating the GPR model (lines 15-17 in Algorithm~\ref{algo:GP_UCB_POMDP}) or random sampling, each experiment is repeated 10 times. All experiments were implemented in Python and solved on a Dell Desktop XPS 8940 with an Intel Core i5-11400 Processor and 40.0 GB RAM to ensure a fair comparison.
\begin{table}[h!]
\centering
\caption{Problem size variables for the test problems}
\label{ExampleProperty}
\begin{tabular}{c|cccc}
\hline
 & $|\cal S|$ & $|\cal A|$ & $|\cal O|$ & $T$ \\ \hline
ChengD51 & 5 & 3 & 3 & 10, 15, 20, 40 \\
Network & 7 & 4 & 2 & 10, 15, 20, 40 \\
Query & 27 & 3 & 3 & 10, 15, 20, 40 \\
Hallway & 60 & 5 & 21 & 10, 15, 20, 40 \\
Aloha.30 & 90 & 29 & 3  & 10, 15, 20, 40 \\ \hline
\end{tabular}
\end{table}
Experiments are terminated based on two stopping criteria: either exceeding 3000 seconds or reaching a target gap of $g_a$. The target gap $g_a$ is dynamically determined after each iteration based on the magnitude of the upper bounds of the value function for the belief point $\bm{b}_0$ at stage $t=0$. Specifically, $g_a$ is computed as
$g_a = \frac{\Lambda \left ( \overline{V_0}(\bm{b}_0) \right )}{10^\rho}~,$
where $\overline{V_0}(\bm{b}_0)$ is the upper bound of the value function for the initial belief point $\bm{b}_0$ at stage $t=0$ after the latest iteration, and $\Lambda(\overline{V_0}(\bm{b}_0))$ rounds $\overline{V_0}(\bm{b}_0)$ up to the nearest power of 10. For example, if $\overline{V_0}(\bm{b}_0) = 281$, then $\Lambda(281) = 1000$, and if $\overline{V_0}(\bm{b}_0) = 64$, then $\Lambda(64) = 100$.
This way, the precision parameter~$\rho$ directly controls the stringency of~$g_a$. Larger values of $\rho$ yield smaller gaps, resulting in more precise approximations. We set the precision level, $\rho$, as $5$ in our experiments. We set the uncertainty level, $\eta$, of the upper confidence bound to $1$, and the accuracy level, $\nu$, of the ALD condition to $10^{-5}$. We set the initial iteration phase in line 11 of Algorithm~\ref{algo:GP_UCB_POMDP} at 5 iterations and perform periodic checks every 5 iterations in our experiments.

To compare and contrast the performance of the different algorithms, we focus on two primary metrics: (i) the lower bound of the value function for the initial belief point $\bm{b_0}$ at stage $t=0$, and (ii) the gap between the upper and lower bounds for $\bm{b_0}$. For experiments involving random selection processes, such as the random selection of belief points in GPR model updates in the GP-UCB approach, both metrics of lower bound and gap are averaged over 10 runs to account for variability and to ensure reliable comparisons. 

In our study, we evaluate the performance of our GP-UCB approach in comparison to the sawtooth projection method under three different belief point sampling strategies: (1) max-gap sampling based on the work from~\cite{walraven2019point}, (2) random sampling proposed by~\cite{spaan2005perseus}, and (3) fixed-grid sampling from~\cite{lovejoy1991computationally}. These strategies are chosen to highlight different aspects of belief point selection, with max-gap focusing on regions of high uncertainty in the value function, random sampling providing a baseline of unbiased coverage, and fixed-grid offering a structured exploration of the belief space. Combining each sampling strategy with the two approximation methods results in a total of six experimental settings, enabling a comprehensive analysis of their relative performance under varying sampling approaches. Among these six settings, fixed-grid sampling is fully deterministic for both the GP-UCB and sawtooth approaches, while max-gap sampling is also deterministic when used with the sawtooth approach.

\section{Performance Results}
\label{sec:perfresults}

\begin{figure}[htp]
    \centering
    \includegraphics[width=1\textwidth]{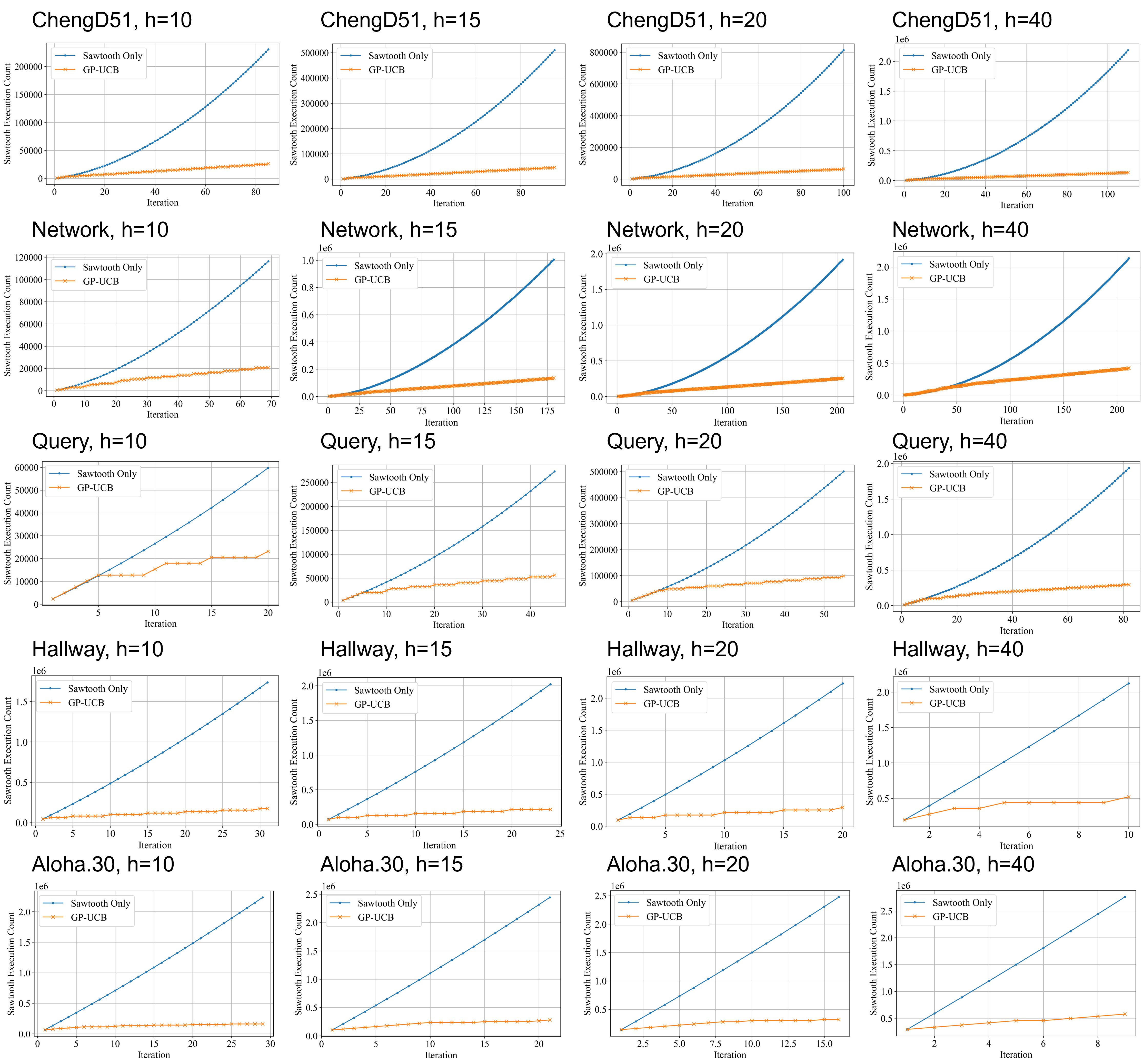}
    \caption{Comparison of the number of sawtooth executions over iterations between algorithms with pure sawtooth and with GP-UCB.}
    \label{fig:SawtoothCounts}
\end{figure}

A key advantage of the GP-UCB approach is its ability to significantly reduce the computational time required for sawtooth projections. Figure~\ref{fig:SawtoothCounts} compares the growth in the number of sawtooth executions over iterations for the max-gap-sawtooth and max-gap-GPUCB approaches, for all test problems. The results indicate that although both methods show an increase in sawtooth executions as iterations progress, the growth is substantially slower for max-gap-GPUCB. Using GP-UCB to estimate the upper bounds for belief points minimizes the need for frequent and computationally expensive sawtooth calculations. This efficiency allows the GP-UCB approach to achieve smaller gaps within the same computational time, or reach the target gap significantly faster than max-gap-sawtooth.

Table~\ref{tab:experimentresults_fiveexamplesh40} presents the experimental results for lower bounds and gaps at $h=40$, based on the predefined stopping criteria. Full results for all problems under all tested planning horizons (i.e., 10, 15, 20, and 40 periods) are provided in Tables~\ref{tab:experimentresults_ChengD51} to~\ref{tab:experimentresults_Aloha30}. Bolded numbers indicate the largest lower bounds (LB) and smallest gaps (Gap), highlighting the best-performing methods for these metrics. The values shown in parentheses represent the standard deviations reflecting the variability observed across multiple runs, where applicable. In addition to reporting the average gap, we also report on the worst-case gap observed across the ten runs to provide additional information on the robustness and consistency of the methods under varying conditions.

Upon comparing the results of GP-UCB and pure sawtooth projection under the max-gap belief expansion strategy, it is evident that GP-UCB consistently achieves superior or comparable performance. Specifically, in terms of the lower bound, GP-UCB reaches the same or higher values as the sawtooth method in all cases by the termination of the experiment runs. Looking at the gaps, GP-UCB outperforms the sawtooth method in all examples, with the exception of the network problem at $h=10$ (full results are provided in Appendix~\ref{sec:fullcompresults}). In this particular case, GP-UCB requires only $6$ additional seconds to reach the target gap, demonstrating no significant disadvantage. For the remaining cases, GP-UCB either reaches the target gap faster or achieves a smaller gap when terminated at the 3000-second limit. Even when evaluating the worst-case gaps, the conclusions remain consistent with the observations for the average gaps across ten runs, further highlighting the efficiency of GP-UCB. 

When evaluating the performance of sawtooth and GP-UCB under the random sampling strategy, we again observe consistent trends. For all test problems, GP-UCB achieves larger lower bounds and smaller gaps compared to the sawtooth method. This highlights the advantage of GP-UCB in improving the value function approximation even when belief points are sampled without a specific focus on regions of high uncertainty.

Finally, when comparing the GP-UCB and sawtooth projection methods under the fixed-grid belief sampling strategy, we observe that the fixed-grid strategy is unable to solve large-scale problems such as query, hallway, and aloha.30. This limitation arises due to the exponential growth in the number of grid points, leading to excessive computational and memory requirements. For the remaining two smaller problems, ChengD51 and network, the GP-UCB method consistently achieves the same lower bounds as the sawtooth method. However, GP-UCB demonstrates its superiority by achieving smaller gaps in all cases. 

The above results demonstrate that the choice of belief expansion strategy significantly impacts performance. Among the three strategies evaluated (i.e.,max-gap, random sampling, and fixed-grid sampling), fixed-grid strategy consistently unperformed and was unable to solve large-scale problems such as query, hallway, and aloha.30 due to its poor scalability to high-dimensional belief spaces. Max-gap sampling generally achieves the smallest gaps, emphasizing its effectiveness in reducing uncertainty in the belief space. However, for large-scale problems like aloha.30, random sampling achieves larger lower bounds than max-gap sampling, as it favors exploration of the belief space.

\begin{table}
\centering
\footnotesize
\renewcommand{\arraystretch}{1.0} % Adjust for vertical spacing
\setlength{\tabcolsep}{4pt} % Adjust column separation if desired
\caption{Experiment results for the five problems at $h=40$. Bolded numbers indicate the largest upper bounds and smallest gaps, values in parentheses are standard deviations.}
\label{tab:experimentresults_fiveexamplesh40}
\begin{tabular}{@{\hspace{5pt}}c@{\hspace{5pt}} @{\hspace{5pt}}l@{\hspace{5pt}} @{\hspace{5pt}}l@{\hspace{5pt}}
                @{\hspace{5pt}}l@{\hspace{5pt}} @{\hspace{5pt}}l@{\hspace{5pt}} @{\hspace{5pt}}l@{\hspace{5pt}}
                @{\hspace{5pt}}l@{\hspace{5pt}} @{\hspace{5pt}}l@{\hspace{5pt}}}
\hline
\multicolumn{2}{@{\hspace{5pt}}c@{\hspace{5pt}}}{} & \makecell[l]{Max-gap\\GP-UCB} & \makecell[l]{Max-gap\\sawtooth} & \makecell[l]{Random\\GP-UCB} & \makecell[l]{Random\\sawtooth} & \makecell[l]{Fixed-grid\\GP-UCB} & \makecell[l]{Fixed-grid\\sawtooth} \\[2pt]\hline

\multirow{5}{*}{ChengD51} & LB 
    & \textbf{261.713}\textsuperscript{(0.000)} 
    & \textbf{261.713} 
    & \textbf{261.713}\textsuperscript{(0.000)}
    & \textbf{261.713}\textsuperscript{(0.000)} 
    & 261.711 
    & 261.711 \\[2pt]
& Gap 
    & \textbf{0.129}\textsuperscript{(0.002)} 
    & 0.170 
    & 2.087\textsuperscript{(0.058)} 
    & 2.287\textsuperscript{(0.092)} 
    & 73.452 
    & 75.449 \\[2pt]
& Time (s) 
    & 3022.0\textsuperscript{(9.3)} 
    & 3028.8 
    & 3069.6\textsuperscript{(26.9)} 
    & 3019.3\textsuperscript{(22.2)} 
    & 2.3 
    & 0.5 \\[2pt]\cdashline{2-8}
& \makecell[c]{Worst Gap}  
    &  \textbf{0.131}
    &  0.170
    &  2.170
    &  2.447
    &  73.452
    &  75.449\\[2pt]\hline

\multirow{5}{*}{Network} & LB 
    & \textbf{592.274}\textsuperscript{(0.003)} 
    & \textbf{592.274}
    & 591.878\textsuperscript{(0.150)} 
    & 591.719\textsuperscript{(0.229)} 
    & 582.704
    & 582.704\\[2pt]
& Gap 
    & \textbf{1.453}\textsuperscript{(0.018)} 
    & 2.419
    & 42.844\textsuperscript{(1.844)} 
    & 50.222\textsuperscript{(1.503)} 
    & 543.012 
    & 660.788\\[2pt]
& Time (s) 
    & 3030.6\textsuperscript{(36.2)} 
    & 3020.2
    & 3015.0\textsuperscript{(15.7)} 
    & 3020.2\textsuperscript{(13.1)} 
    & 3.2 
    & 1.5\\[2pt]\cdashline{2-8}
& \makecell[c]{Worst Gap}  
    &  \textbf{1.488}
    &  2.419
    &  45.825
    &  53.214
    &  543.012
    &  660.788\\[2pt]\hline

\multirow{5}{*}{Query} & LB 
    & \textbf{120.327}\textsuperscript{(0.000)} 
    & \textbf{120.327}
    & \textbf{120.327}\textsuperscript{(0.000)}
    & \textbf{120.327}\textsuperscript{(0.000)} 
    & NA & NA \\[2pt]
& Gap 
    & 0.017\textsuperscript{(0.001)} 
    & 0.026
    & \textbf{0.010}\textsuperscript{(0.001)}
    & 0.037\textsuperscript{(0.000)} 
    & NA & NA \\[2pt]
& Time (s) 
    & 3004.4\textsuperscript{(5.6)} 
    & 3049.5
    & 3018.8\textsuperscript{(32.9)}
    & 3027.4\textsuperscript{(30.6)} 
    & NA & NA \\[2pt]\cdashline{2-8}
& \makecell[c]{Worst Gap}  
    &  0.018
    &  0.026
    &  \textbf{0.010}
    &  0.037
    &  NA
    &  NA \\[2pt]\hline

\multirow{5}{*}{Hallway} & LB 
    & \textbf{1.930}\textsuperscript{(0.000)} 
    & 1.912
    & 1.612\textsuperscript{(0.067)} 
    & 1.295\textsuperscript{(0.061)} 
    & NA & NA \\[2pt]
& Gap 
    & \textbf{1.220}\textsuperscript{(0.000)} 
    & 1.234
    & 1.451\textsuperscript{(0.070)} 
    & 1.510\textsuperscript{(0.061)} 
    & NA & NA \\[2pt]
& Time (s) 
    & 2871.8\textsuperscript{(18.9)} 
    & 3030.7
    & 3092.7\textsuperscript{(66.6)}  
    & 3069.2\textsuperscript{(21.0)} 
    & NA & NA \\[2pt]\cdashline{2-8}
& \makecell[c]{Worst Gap}  
    &  \textbf{1.220}
    &  1.234
    &  1.571
    &  1.638
    &  NA
    &  NA \\[2pt]\hline

\multirow{5}{*}{Aloha.30} & LB 
    & 282.129\textsuperscript{(0.000)} 
    & 282.039
    & \textbf{283.338}\textsuperscript{(0.249)}
    & 282.329\textsuperscript{(0.577)} 
    & NA & NA \\[2pt]
& Gap 
    & \textbf{10.428}\textsuperscript{(0.003)} 
    & 12.507
    & 11.164\textsuperscript{(0.632)}
    & 12.480\textsuperscript{(0.932)} 
    & NA & NA \\[2pt]
& Time (s) 
    & 2956.0\textsuperscript{(27.2)} 
    & 3086.8
    & 3084.9\textsuperscript{(45.8)}
    & 3163.7\textsuperscript{(39.4)} 
    & NA & NA \\[2pt]\cdashline{2-8}
& \makecell[c]{Worst Gap}  
    &  \textbf{10.435}
    &  12.507
    &  12.113
    &  14.570
    &  NA
    &  NA \\[2pt]\hline
\end{tabular}
\end{table}

To further evaluate the efficiency of GP-UCB compared to the sawtooth projection method, we analyze the time required for GP-UCB to achieve the same (or smaller) gap as sawtooth at its terminating iteration under the max-gap belief expansion strategy. The results, summarized in Table~\ref{tab:time_gap_problems}, demonstrate significant time savings in most problems and planning horizons. For smaller problems like ChengD51 and network, GP-UCB achieves the sawtooth gaps in 30\%-60\% less time. For larger problems such as query, hallway, and aloha.30, the time savings are even more significant, reaching up to 99.7\% for query and over 70\% for aloha.30 at longer horizon lengths. While GP-UCB generally provides greater time savings at longer horizons, some variations exist across problems, such as Hallway. These differences likely stem from problem-specific structures and how efficiently GP-UCB generalizes across horizons. To capture these trends and exceptions, we present results for all planning horizons in Table~\ref{tab:time_gap_problems}.

\begin{table}
    \centering
    \footnotesize
    \caption{Time for GP-UCB to reach the same or smaller gap as sawtooth achieved at the terminating iteration for different problems using max-gap method for belief expansion}
    \label{tab:time_gap_problems}
    \renewcommand{\arraystretch}{1} % Adjust row height
    \begin{tabular}{l c c c c c c}
        \hline
        \multirow{2}{*}{Problem} & \multirow{2}{*}{Planning Horizon} 
        & \multicolumn{2}{c}{Sawtooth} & \multicolumn{2}{c}{GP-UCB} & \multirow{2}{*}{Time Reduction (\%)}\\
        \cline{3-6}
        & & Gap & Time (s) & Gap & Time (s) & \\
        \hline
        \multirow{4}{*}{ChengD51} & h=10 & 0.005 & 3022.4 & 0.005 & 1860.5 & 38.4\%\\
                                   & h=15 & 0.021 & 3010.9  & 0.021 & 1873.5 & 37.8\%\\
                                   & h=20 & 0.045 & 3032.9  & 0.044 & 1685.1 & 44.4\%\\
                                   & h=40 & 0.170 & 3028.8  & 0.167 & 1909.0 & 37.0\%\\
        \hline
        \multirow{4}{*}{Network} & h=10 & 0.010 & 52.9 & 0.010 & 58.8 & -11.2\%\\
                                   & h=15 & 0.039 & 3011.3  & 0.038 & 2060.9 & 31.6\%\\
                                   & h=20 & 0.346 & 3029.9  & 0.342 & 1183.4 & 60.9\%\\
                                   & h=40 & 2.419 & 3020.2  & 2.340 & 1369.1 & 54.7\%\\
        \hline
        \multirow{4}{*}{Query} & h=10 & 0.001 & 80.7 & 0.001 & 5.7 & 92.9\%\\
                                   & h=15 & 0.002 & 3000.8  & 0.001 & 8.6 & 99.7\%\\
                                   & h=20 & 0.004 & 3025.6  & 0.003 & 11.1 & 99.6\%\\
                                   & h=40 & 0.026 & 3049.5  & 0.026 & 33.5 & 98.9\%\\
        \hline
        \multirow{4}{*}{Hallway} & h=10 & 0.168 & 3034.5 & 0.166 & 2759.3 & 9.1\%\\
                                   & h=15 & 0.360 & 3082.8  & 0.358 & 1458.9 & 52.7\%\\
                                   & h=20 & 0.524 & 2883.4  & 0.524 & 2337.9 & 18.9\%\\
                                   & h=40 & 1.234 & 3030.7  & 1.222 & 2922.8 & 3.6\%\\
        \hline
        \multirow{4}{*}{Aloha.30} & h=10 & 0.736 & 2962.7 & 0.735 & 2728.3 & 7.9\%\\
                                   & h=15 & 2.004 & 2912.5  & 1.978 & 776.9 & 73.3\%\\
                                   & h=20 & 4.003 & 3070.2  & 3.952 & 937.7 & 69.5\%\\
                                   & h=40 & 12.507 & 3086.8  & 12.419 & 740.2 & 76.0\%\\
        \hline
    \end{tabular}
\end{table}

\begin{figure}[h]
  \centering
  \includegraphics[width=\linewidth]{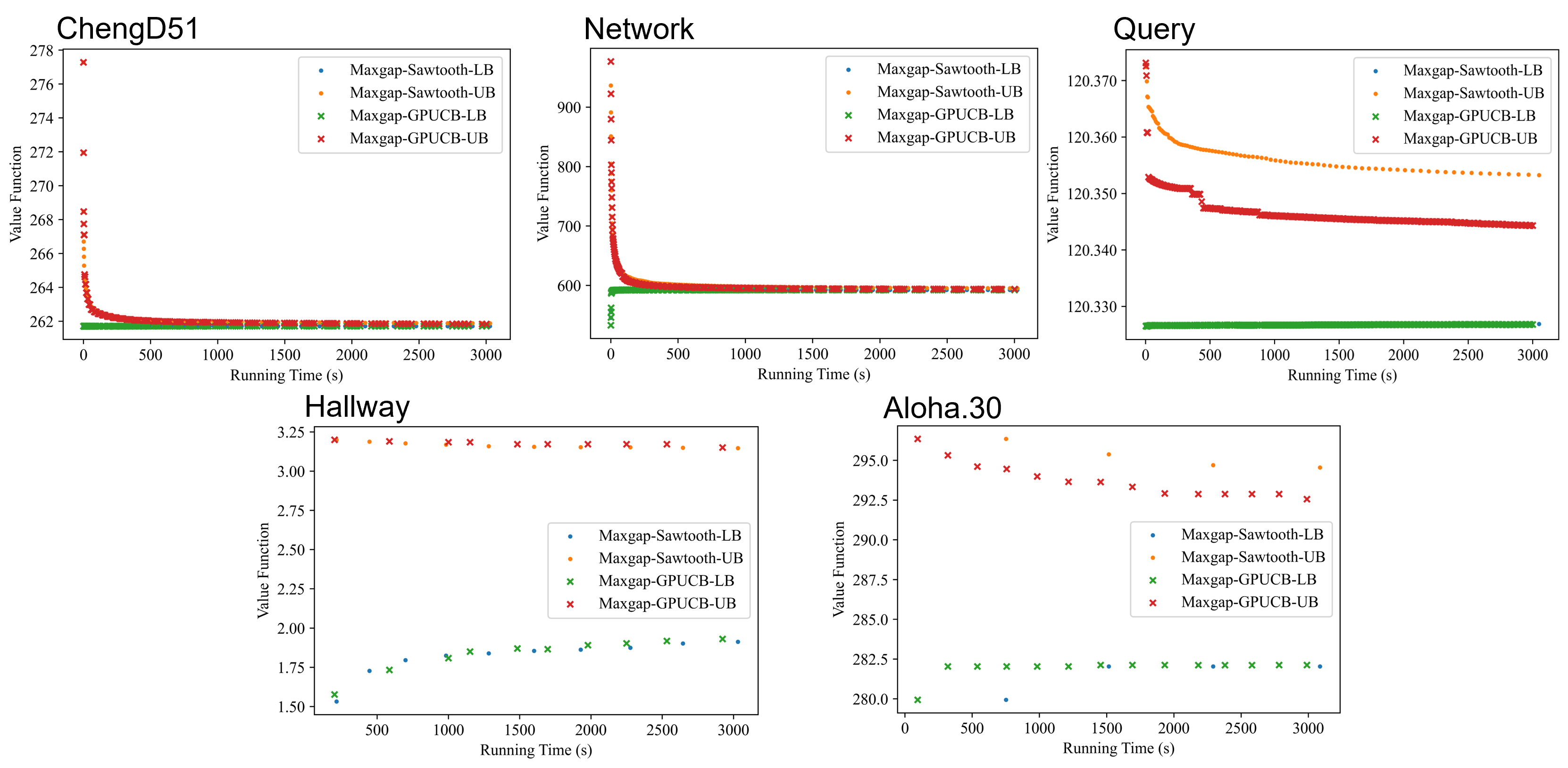}
  \caption{Upper and lower bounds of the value function for the five examples at $h=40$. Dots represent upper and lower bounds from FiVI, while crosses indicate those from AUCB, color-coded as shown in the legend.}
  \label{fig:ValueFunction_h40}
\end{figure}

To provide a visual representation of the experiment results, we plot the lower bounds and upper bounds of the value function during the execution of the max-gap-GP-UCB and max-gap-sawtooth methods at $h=40$ in Figure~\ref{fig:ValueFunction_h40}. Full results are shown in Figure~\ref{fig:ValueFunction}. Consistent with the findings in Tables~\ref{tab:experimentresults_fiveexamplesh40}, the GP-UCB approach demonstrates strong performance in all test problems, achieving more favorable lower and upper  bounds in most cases. For problems such as query and aloha.30, GP-UCB achieves significantly smaller upper bounds from the very beginning of the experiment runs, demonstrating its efficiency in approximating the value function in these cases. In contrast, for problems like the hallway problem, the reduction in upper bounds is less significant, indicating that certain problem structures may inherently limit the performance of GP-UCB.

One limitation of our approach is the need to select the kernel function. In these experiments, we used an exponential kernel for the GP-UCB approach. While this choice works well for most problems, larger problems like hallway and aloha.30 may benefit from alternative kernels such as the Matérn kernel, which can better generalize across larger belief spaces and avoid overfitting to the small initial data sets.

\section{Conclusions}
\label{sec:Conclusion}

In this paper, we present GP-UCB, a novel approach to accelerate the upper bound estimation in point-based value iteration to solve finite-horizon POMDPs. By utilizing Gaussian Process Regression to approximate the upper bound convex hull across the belief space, our method reduces computational complexity without sacrificing solution quality. The approach selectively trains on the most informative subset of belief points, significantly enhancing solution efficiency and scalability. This innovation avoids redundant upper bound updates and provides theoretical guarantees for the convergence of the proposed GP-UCB to the true upper bound convex hull. Through extensive experiments on benchmark problems, we demonstrated that GP-UCB consistently accelerates convergence and improves scalability, reducing computation time by 30\%-60\% on smaller problems and up to 99.7\% on larger ones while maintaining the same upper bound gaps as the pure sawtooth projection method. These results underscore its potential as a powerful tool for solving practical large-scale POMDPs in finite-horizon settings.

In the proposed approach, the GPR model uses training data from sawtooth projections, which are computationally efficient for estimating upper bounds. However, a trade-off exists between computation time and data quality. Generating higher-quality data results in better approximations for upper bounds, but may take longer. Future research will focus on finding a way to balance this trade-off to improve GPR model performance. Another promising direction is optimizing the selection of kernels for GPR, as the choice of kernel directly impacts the quality of the approximation and computational efficiency. Investigating adaptive kernel selection strategies or using domain-specific kernels could enhance the flexibility and accuracy of the GPR model.

% Acknowledgments here
%\ACKNOWLEDGMENT{%
% Enter the text of acknowledgments here
%}% Leave this (end of acknowledgment)

% Appendix here
% Options are (1) APPENDIX (with or without general title) or 
%             (2) APPENDICES (if it has more than one unrelated sections)
% Outcomment the appropriate case if necessary
%
% \begin{APPENDIX}{<Title of the Appendix>}
% \end{APPENDIX}
%
%   or 
%

\newpage
\appendix

\section{Background on POMDP Problems and Solution Algorithms}
\label{sec:background}

Like MDPs, POMDPs identify control actions based on the current state of the system, but since the system state is only partially observable, decisions are made based on a \textit{belief state}—a probability distribution over the possible system states that represents the current estimate of the system’s true state. This belief state is updated through systemic observations that provide probabilistic information about the system state at any given time. In addition to the usual MDP assumptions regarding action- and state-dependent transition probabilities, POMDP formulations assume a \textit{likelihood function}, which gives the probability of observing a particular outcome for each system state. This allows for the computation of the next belief state as a function of the current belief state, the action taken, and the observation obtained~\citep{littman2009tutorial}. Since the belief state is a continuous vector over the probability distribution of the system states, the belief space becomes continuous, making the solution of POMDPs much more difficult~\citep{papadimitriou1987complexity, kaelbling1998planning, madani1999undecidability}.

The structural properties of the \textit{value function} in POMDPs, for both finite- and infinite-horizon problems, are well understood. In both cases, the value function, which maps a belief state to the expected total reward, is \textit{piecewise-linear and convex} in the belief space~\citep{sondik1971optimal, smallwood1973optimal}. This means that the value function can be represented as the maximum of a set of linear functions over the belief space, each associated with a specific action. These linear functions are known as \textit{alpha-vectors}, where each alpha-vector corresponds to a specific action and encodes the expected future rewards for that action given the current belief state~\citep{sondik1971optimal}. For any belief state, the value is determined by selecting the alpha-vector that maximizes the expected reward. In the finite-horizon case, the set of alpha-vectors is updated stage by stage using backward induction to compute the value function at each time step~\citep{smallwood1973optimal}. In infinite-horizon problems, the value function is stationary, and the alpha-vectors are iteratively refined until a fixed point is reached~\citep{sondik1978optimal, monahan1982survey, lovejoy1991survey}. However, the number of belief states and corresponding alpha-vectors grows \textit{exponentially} with the size of the state and observation spaces, leading to the well-known ``curse of dimensionality.'' This exponential growth makes the computation of exact solutions infeasible for large-scale problems, even when pruning techniques are employed to reduce the number of alpha-vectors~\citep{sondik1971optimal,monahan1982survey,cheng1988algorithms,littman1996algorithms,zhang1996planning,cassandra1997incremental}.

\subsection{Point-Based Value Iteration}
After the structural properties of POMDP value functions were established, and the challenges of solving them exactly for large problems became evident, researchers developed various approximate methods. Among the most prominent of these is the Point-Based Value Iteration (PBVI) algorithm, introduced by \citet{pineau2003pbvi}. PBVI was a major breakthrough because it significantly improved the scalability of POMDP solvers by focusing on a subset of belief points, rather than attempting to compute the value function over the entire continuous belief space. 

The key insight behind PBVI is that not all belief states are equally important for decision making. Instead of considering the entire belief space, PBVI selects a representative set of belief points and performs computations on these points. By concentrating computational effort on a smaller set of critical belief states, PBVI drastically reduces the computational load while still providing a good approximation of the value function. This approach is particularly effective, as PBVI uses the structure of POMDPs where a limited set of belief points - especially those likely to be encountered under the optimal policy - can be sufficient to approximate the value function well enough for effective decision making~\citep{pineau2003pbvi, smith2005point}. This makes PBVI one of the most widely used algorithms to solve POMDPs, particularly for large and complex problems~\citep{shani2013survey}.  

The main components of PBVI algorithms can be summarized as follows:
\begin{enumerate}
\item {\bf Select Initial Belief Points:} A small set of representative belief points is initialized, typically using random sampling or simulations.
\item {\bf Backup Operation:} The lower and upper bounds on the value function at each belief point are updated by considering the immediate reward and expected future rewards based on potential actions and observations.
\item {\bf Expand Belief Set:} New belief points are added as needed to improve the approximation across the belief space.
\item {\bf Convergence:} The process repeats until the difference between the upper and lower bounds falls below a set threshold, indicating convergence.
\end{enumerate}

Since its introduction in 2003, various improvements have been made to each of the four key components of PBVI. For belief point selection, techniques such as heuristic sampling~\citep{pineau2003pbvi} and more structured approaches such as reachability analysis~\citep{kurniawati2008sarsop} have improved the initial set of belief points by focusing on those most relevant to the policy. The backup operation has seen optimizations with randomized belief backups, as in the Perseus algorithm~\citep{spaan2005perseus}, which selectively updates belief points to reduce computational cost. In terms of belief set expansion, advanced methods such as adaptive belief selection~\citep{shani2013survey} and forward search techniques have improved the algorithm's ability to refine the belief set efficiently.

For convergence, several approaches have been employed to accelerate the process. One such strategy is incremental pruning~\citep{cassandra1997incremental}, which reduces the complexity of managing alpha-vectors by eliminating irrelevant vectors during each iteration. Another is heuristic search, which focuses computational resources on the most promising regions of the belief space, improving both convergence speed and efficiency~\citep{smith2004heuristic}. Additionally, bounded policy updates as implemented in Heuristic Search Value Iteration (HSVI)~\citep{smith2005point} ensure that the value function and policy are refined efficiently while maintaining convergence guarantees.

\subsection{Upper and Lower Bounds}
One of the most significant innovations for convergence and belief point selection is the introduction of upper and lower bounds on the value function. These bounds serve multiple purposes: they guide the selection of critical belief points by focusing on regions where the difference between the upper and lower bounds is largest, and they provide stopping criteria for the algorithm when the bounds converge sufficiently~\citep{smith2005point, poupart2011closing}. By integrating bounds into PBVI, convergence is not only accelerated, but the algorithm also gains the ability to focus more intelligently on belief points where improvement is needed most.

In PBVI algorithms, calculating the lower bounds is relatively straightforward and computationally inexpensive. For each belief point, the best \textit{alpha-vector} (the one that yields the highest expected reward) is chosen, ensuring that the lower bound represents a conservative estimate of the value function. This method provides a reliable, minimal reward estimate without adding significant computational overhead. However, calculating upper bounds has traditionally been more complex and costly. Exact methods for the upper bound calculation involve projecting the belief point onto the convex hull of the belief/upper bound pairs. The tightest upper bound is obtained by minimizing the projection using linear programming methods. Although this provides the strongest possible upper bound, it is computationally expensive~\citep{monahan1982survey, cheng1988algorithms}. The QMDP approach, which uses $Q$-values (i.e., the expected reward of taking an action in a state and following an optimal policy), simplifies the problem by assuming full observability after a single step, leading to an optimistic upper bound~\citep{littman1995witness}. 

For such cases, approximate methods like those used in HSVI and Successive Approximations of the Reachable Space under Optimal Policies (SARSOP) algorithms offer more computational efficiency by deriving upper bounds from optimistic assumptions. In HSVI, the upper bound is updated by a heuristic search, exploring actions and observations optimistically, guiding the algorithm toward promising belief points~\citep{smith2005point}. SARSOP focuses on reachable belief spaces to update upper bounds only for belief points likely to be encountered under an optimal policy, further reducing computational costs~\citep{kurniawati2008sarsop}. These heuristic-based methods are particularly useful in infinite-horizon POMDPs, where the value function is stationary, meaning upper bounds do not need frequent recalculation. 

In finite-horizon POMDPs, the situation becomes more complex because the value function changes dynamically over time~\citep{walraven2019point}. The upper bound must be recalculated at each time step, reflecting the evolving decision horizon. Calculating an exact upper bound for finite-horizon problems involves exploring all possible future outcomes optimistically, considering the best-case scenario at each stage, which leads to a high computational cost~\citep{walraven2019point, smith2005point}. This process requires evaluating future rewards under all potential actions and observations for each belief point, further increasing the complexity.

Although PBVI has been effective for infinite-horizon problems, where a stationary policy optimizes the value function indefinitely, significant adjustments are necessary to apply it to finite-horizon problems. In infinite-horizon settings, the value function remains stationary, allowing algorithms such as PBVI to refine the value function iteratively without accounting for the remaining time steps~\citep{pineau2003pbvi}. However, in finite-horizon problems, the value function changes since the number of remaining decisions affects both immediate and future rewards. As the horizon shortens, the value function increasingly reflects immediate rather than future rewards. This requires the algorithm to adapt at each step to account for the evolving optimal policy~\citep{walraven2019point}. Efficient mechanisms, particularly for accurate upper and lower bounds, are crucial to ensuring that belief point selection and value function updates remain computationally feasible.

To address these challenges, \citet{walraven2019point} implemented efficient methods for handling upper and lower bounds in finite-horizon problems. By estimating the maximum possible value for each belief point and focusing on the widest gap between the upper and lower bounds, the algorithm prioritizes regions where the value function is the least accurate. This accelerates convergence and reduces computation time, while upper bounds serve as a stopping criterion, terminating the algorithm when convergence is reached, ensuring a near-optimal policy without excessive refinement.

\subsection{Sawtooth Projection to Obtain Upper Bounds}

Instead of using exact upper bounds,  \citet{hauskrecht2000value} introduced the so-called {\it sawtooth projections} to dynamically adjust the upper bounds at each step of the horizon. This approximation leverages the convex structure of the value function to reduce computational complexity while maintaining a reasonable level of approximation accuracy. \citet{poupart2011closing} reported that sawtooth projections can significantly reduce computational costs while preserving a high degree of solution quality. Specifically, in the worst-case scenario, the computational complexity for the approximation decreases to $O(|\cal B||\cal S|)$, where $\cal{B}$ is a sampled belief set in the $|\cal S|$-dimensional space with $|\cal B|>|\cal S|$. 

The effectiveness of sawtooth projection in finite-horizon settings was recently demonstrated by~\citet{walraven2019point} in their Finite-horizon Value Iteration (FiVI) algorithm. FiVI outperforms both the original PBVI~\citep{pineau2003pbvi} and the Heuristic Search Value Iteration (HSVI)~\citep{smith2005point} in terms of speed and scalability, especially in large finite-horizon problems. By focusing computational effort where it is most needed, the sawtooth projection allows FiVI to achieve faster convergence without sacrificing solution quality. Furthermore, their numerical results showed that FiVI scales better with both the horizon length and the complexity of belief spaces. This makes it a more versatile option for large-scale applications, and less sensitive to the size of state and action spaces~\citep{walraven2019point}.

The primary computational complexity of the sawtooth projection stems from the need to repeatedly update the upper bound each time a new belief point is added to the belief set. As this set expands, the upper bound approximation of the value function improves, narrowing the gap with the lower bound. In the initial stages of belief expansion, adding new points significantly improves the upper bound. However, as more belief points are added, the marginal improvement becomes minimal, even though the number of updates ---and thus the computational complexity---increases \citep{hauskrecht2000value}. This leads to upper bound updates that require computations that scale with the product of the cardinality of the belief set and the dimensionality of the belief state, resulting in a computational complexity of $O(|\cal B||\cal S|)$.

Hence, while sawtooth projection offers a computationally efficient alternative to computing the convex hull of upper bounds, its effectiveness remains limited, as the number of sawtooth projection operations grows rapidly with the expanding belief set. As a result, this approach is still insufficient to significantly reduce complexity and solve larger problems.

In this work, we model the convex hull of the upper bounds by fitting a subset of the belief points and their sawtooth projections using GPR. As such, for any new belief point, the upper bound interpolation is directly inferred from the posterior of the GPR predictions, precluding the need to solve any linear programs or sawtooth projection. Furthermore, at every iteration of the upper bound approximation, we fit the GPR only to a subset of the belief/upper bound pairs that are most informative, to reduce the computational cost of covariance matrix inversion. 

Notably, multiple belief points lying on the same alpha-vector are redundant and do not provide additional information about the value function. To avoid such redundant belief points while training the GPR, we argue that the most informative belief points are those that are approximately linearly independent in a high-dimensional Hilbert space \citep{islam2024dynamic, engel2004kernel}. Therefore, each time a new belief point is added, we select a subset of belief points that satisfy the approximate linear dependence (ALD) criterion to train the GPR. This selected subset is usually significantly smaller than the entire belief set, which further reduces computational demands.

Finally, it must be noted that the proposed GP-UCB is only approximately correct, although with a very high probability. Hence, we provide theoretical results on the consistency of the proposed upper bound approximations as the belief set expands.
\section{Description of the tiger problem}
\label{TigerProblem}
To demonstrate the concept of reachable belief points in a POMDP, we use the tiger problem, originally introduced by~\citet{kaelbling1998planning}. In this problem, an agent faces two doors, with a tiger behind one and treasure behind the other, i.e., the unobservable states are ``Tiger Left'' and ``Tiger Right.'' At every stage, the agent can open one of the doors (i.e., actions ``open left" or ``open right'') or listen (action ``listen'') for a tiger growl to gather (noisy) information. In particular, there is an 85\% chance of correctly estimating the location of the tiger after listening. The agent’s actions are associated with specific rewards and costs: opening the correct door yields a reward of $+10$, while opening the door with the tiger results in a penalty of $-100$, and listening incurs a cost of $-1$. After the agent opens a door, the state resets: the location of the tiger is randomized, placing it behind one of the doors with equal probability. Observations are then made in the new state, where both left and right actions lead to either observation (Growl from Left or Growl from Right) with a probability of 0.5, regardless of the actual location of the tiger.  In the finite-horizon version, at the terminating stage, only the immediate reward or penalty from the final action is considered, without further actions or resets.

\section{Further details on sawtooth projection method}
\label{sec:sawtooth}

\setcounter{equation}{0} % Reset equation counter
\renewcommand{\theequation}{C.\arabic{equation}} % Change numbering to B.1, B.2

\renewcommand{\thealgorithm}{C.\arabic{algorithm}} % Numbering as B.1, B.2, etc.
\setcounter{algorithm}{0} % Restart numbering

\renewcommand{\thefigure}{C.\arabic{figure}}
\setcounter{figure}{0}

\citet{hauskrecht2000value} proposed an interpolation approach for efficient approximation of the convex hull, referred to as {\it sawtooth projection}, to reduce computational complexity. Given an arbitrary belief $\bm{b} \in \cal{B}$ and corner beliefs $\{\bm{w}_1,\bm{w}_2,...,\bm{w}_{|\cal{S}|}\}$, the sawtooth projection provides the upper bound interpolation for any new belief point $\bm{b'}$. The corner beliefs are defined as $\bm{w}_s=(b(1),b(2),...b(|\cal{S}|)$, where $b(s)=1$ and $b(s')=0$ for all $s' \neq s$. The best upper bound approximation for a belief $\bm{b}$ is the one that minimizes $\lambda(\bm{b}) f(\bm{b})$, where
\begin{eqnarray}
    \lambda(\bm{b}) &=& \min\{b'(s)/b(s)|s\in \cal S\}~, \text{ and } \\
    f(\bm{b}) &=& \overline{V}(\bm{b}) - \sum_{s\in \cal{S}}b(s)\overline{V}(\bm{w}_s)~,
\end{eqnarray}
Intuitively, this method interprets the value function as a set of downward-pointing pyramids, with their bases at the corners of the belief space and their tips at known belief point upper bounds.

Following the sawtooth projection for an arbitrary belief point $\bm{b'}$ at time stage $t$, the upper bound for belief point, $\bm{b}$ at stage $t$ is obtained as
\begin{equation}
\overline{V_t}(\bm{b}) = \underset{a\in {\cal A}}{\max} \left[ \sum_{s\in {\cal S}} b(s)R(s,a) +  \sum_{o \in {\cal O}} P(o|\bm{b},a)\bar{v}_{t+1}(\bm{b'}) \right]~, \label{eq:upperbound}
\end{equation} 
where $\bar{v}_{t+1}(\bm{b}')$ is the upper bound projection obtained using sawtooth projection for the belief state $\bm{b'}$ resulting from taking action $a$ and observing $o$. Details of the upper bound projection using the sawtooth projection is presented in Algorithm~\ref{Algo:sawtooth}. Intuitively, the upper bound update is the maximum expected value of the sum of immediate reward and the upper bound of the future reward. 

\begin{algorithm}[H]
\small
    \caption{Upperbound approximation with sawtooth projection method}
    \label{Algo:sawtooth}
    \begin{algorithmic}[1] % Enables line numbering
    \State \textbf{Input:} POMDP model, belief point $\bm{b'}$, belief-bound pairs $(\bm{b},\overline{V}(\bm{b}))$ for each belief $\bm{b}$ in the belief set $\cal{B}$
    \State \textbf{Output:} upper bound corresponding to belief point $\bm{b'}$
    \For{$\bm{b} \in \cal{B} \setminus \{\boldsymbol{w}_s | s \in S\}$}
        \State $f(\bm{b}) \leftarrow \overline{V}(\bm{b}) - \sum_{s \in S} b(s)\overline{V}(\boldsymbol{w}_s)$
        \State $\lambda(\bm{b}) \leftarrow \min_s \frac{b'(s)}{b(s)}$
    \EndFor
    \State $\bm{b}^* \leftarrow \arg\min_b \lambda(\bm{b}) f(\bm{b})$
    \State \Return $\lambda(\bm{b}^*) f(\bm{b}^*) + \sum_{s \in S} b'(s) \overline{V}(\boldsymbol{w}_s)$
    \end{algorithmic}
\end{algorithm}

\begin{figure}[htp]
\centering
    \includegraphics[width=1\textwidth]{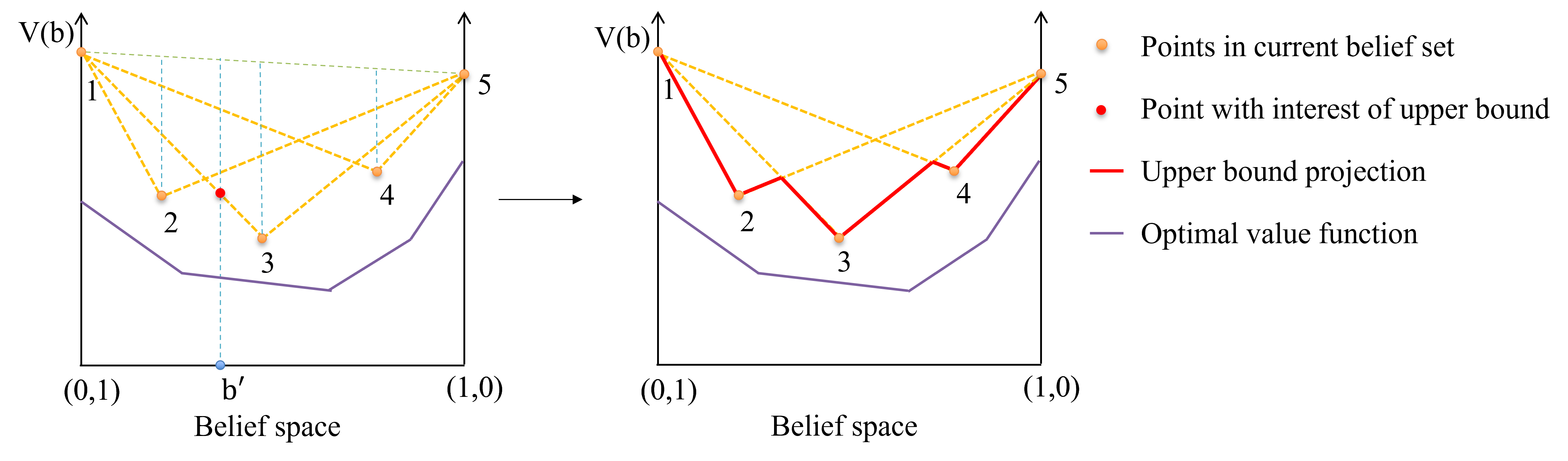}
    \caption{The approximation of upper bound with sawtooth projection.}
    \label{fig:Sawtooth}
\end{figure}

Figure~\ref{fig:Sawtooth} illustrates the sawtooth projection method for a two-state problem. The orange dots labeled 1 to 5 represent belief points with known upper bounds, with dots 1 and 5 being corner beliefs. By connecting belief points 2, 3, and 4 to the corner beliefs, three downward pointing triangles are formed. For a new belief point $\bm{b'}$, the sawtooth method estimates its upper bound by projecting $\bm{b'}$ onto the connecting lines and choosing the projection with the lowest value, which in this case is the line formed by dot~3 and corner belief~1. The estimated upper bound is represented by the red dot on the left of Figure~\ref{fig:Sawtooth}. The sawtooth projection for the entire belief space is shown as the red solid line in the right figure. For comparison, the optimal value function is plotted as the purple line in Figure~\ref{fig:Sawtooth}.

\section{Full computational results}
\label{sec:fullcompresults}

\renewcommand{\thetable}{D.\arabic{table}}
\renewcommand{\thefigure}{D.\arabic{figure}}
\setcounter{table}{0}
\setcounter{figure}{0}

\begin{table}[htp]
\centering
\footnotesize
\renewcommand{\arraystretch}{1.0} % Adjust for vertical spacing
\setlength{\tabcolsep}{4pt} % Adjust column separation if desired
\caption{Experiment results for the ChengD51 problem. The target gap is $0.001$. Bolded numbers indicate the largest upper bounds and smallest gaps, values in parentheses are standard deviations.}
\label{tab:experimentresults_ChengD51}
\begin{tabular}{@{\hspace{5pt}}c@{\hspace{5pt}} @{\hspace{5pt}}l@{\hspace{5pt}} @{\hspace{5pt}}l@{\hspace{5pt}}
                @{\hspace{5pt}}l@{\hspace{5pt}} @{\hspace{5pt}}l@{\hspace{5pt}} @{\hspace{5pt}}l@{\hspace{5pt}}
                @{\hspace{5pt}}l@{\hspace{5pt}} @{\hspace{5pt}}l@{\hspace{5pt}}}
\hline
\multicolumn{2}{@{\hspace{5pt}}c@{\hspace{5pt}}}{} & \makecell[l]{Max-gap\\GP-UCB} & \makecell[l]{Max-gap\\sawtooth} & \makecell[l]{Random\\GP-UCB} & \makecell[l]{Random\\sawtooth} & \makecell[l]{Fixed-grid\\GP-UCB} & \makecell[l]{Fixed-grid\\sawtooth} \\[2pt]\hline
\multirow{5}{*}{h=10} & LB 
    & \textbf{65.246}\textsuperscript{(0.000)} 
    & \textbf{65.246} 
    & \textbf{65.246}\textsuperscript{(0.000)}
    & \textbf{65.246}\textsuperscript{(0.000)} 
    & \textbf{65.246} 
    & \textbf{65.246} \\[2pt]
& Gap
    & \textbf{0.004}\textsuperscript{(0.000)} 
    & 0.005 
    & 0.231\textsuperscript{(0.026)} 
    & 0.263\textsuperscript{(0.020)} 
    & 16.897 
    & 18.896 \\[2pt]
& Time (s) 
    & 3008.2\textsuperscript{(7.0)} 
    & 3022.4 
    & 3013.8\textsuperscript{(8.2)} 
    & 3017.4\textsuperscript{(6.9)} 
    & 0.4 
    & 0.1 \\[2pt]\cdashline{2-8}
& \makecell[c]{Worst Gap}  
    &  \textbf{0.004}
    &  0.005
    &  0.278
    &  0.308
    &  16.897
    &  18.896\\[2pt]\hline

\multirow{5}{*}{h=15} & LB 
    & \textbf{97.990}\textsuperscript{(0.000)} 
    & \textbf{97.990} 
    & \textbf{97.990}\textsuperscript{(0.000)}
    & \textbf{97.990}\textsuperscript{(0.000)} 
    & \textbf{97.990} 
    & \textbf{97.990} \\[2pt]
& Gap 
    & \textbf{0.017}\textsuperscript{(0.000)} 
    & 0.021 
    & 0.488\textsuperscript{(0.014)} 
    & 0.562\textsuperscript{(0.049)} 
    & 26.322 
    & 28.322 \\[2pt]
& Time (s) 
    & 3020.0\textsuperscript{(17.6)} 
    & 3010.9 
    & 3014.7\textsuperscript{(10.1)} 
    & 3018.2\textsuperscript{(9.4)} 
    & 0.7 
    & 0.2 \\[2pt]\cdashline{2-8}
& \makecell[c]{Worst Gap}  
    &  \textbf{0.017}
    &  0.021
    &  0.502
    &  0.650
    &  26.322
    &  28.322\\[2pt]\hline

\multirow{5}{*}{h=20} & LB 
    & \textbf{130.735}\textsuperscript{(0.000)} 
    & \textbf{130.735} 
    & \textbf{130.735}\textsuperscript{(0.000)}
    & \textbf{130.735}\textsuperscript{(0.000)} 
    & 130.734 
    & 130.734 \\[2pt]
& Gap 
    & \textbf{0.034}\textsuperscript{(0.001)} 
    & 0.045 
    & 0.758\textsuperscript{(0.051)}
    & 0.850\textsuperscript{(0.039)} 
    & 35.748 
    & 37.747 \\[2pt]
& Time (s) 
    & 3019.9\textsuperscript{(21.1)} 
    & 3032.9 
    & 3047.9\textsuperscript{(17.8)}
    & 3033.9\textsuperscript{(8.5)} 
    & 1.1 
    & 0.2 \\[2pt]\cdashline{2-8}
& \makecell[c]{Worst Gap}  
    &  \textbf{0.035}
    &  0.045
    &  0.842
    &  0.913
    &  35.748
    &  37.747\\[2pt]\hline

\multirow{5}{*}{h=40} & LB 
    & \textbf{261.713}\textsuperscript{(0.000)} 
    & \textbf{261.713} 
    & \textbf{261.713}\textsuperscript{(0.000)}
    & \textbf{261.713}\textsuperscript{(0.000)} 
    & 261.711 
    & 261.711 \\[2pt]
& Gap 
    & \textbf{0.129}\textsuperscript{(0.002)} 
    & 0.170 
    & 2.087\textsuperscript{(0.058)} 
    & 2.287\textsuperscript{(0.092)} 
    & 73.452 
    & 75.449 \\[2pt]
& Time (s) 
    & 3022.0\textsuperscript{(9.3)} 
    & 3028.8 
    & 3069.6\textsuperscript{(26.9)} 
    & 3019.3\textsuperscript{(22.2)} 
    & 2.3 
    & 0.5 \\[2pt]\cdashline{2-8}
& \makecell[c]{Worst Gap}  
    &  \textbf{0.131}
    &  0.170
    &  2.170
    &  2.447
    &  73.452
    &  75.449\\[2pt]\hline

\end{tabular}
\end{table}

\begin{table}[htp]
\centering
\footnotesize
\renewcommand{\arraystretch}{1.0} % Adjust this to reduce or increase vertical spacing
\setlength{\tabcolsep}{4pt} % Adjust column separation if desired
\caption{Experiment results for the network problem. The target gap is $0.01$. Bolded numbers indicate the largest upper bounds and smallest gaps, values in parentheses are standard deviations.}
\label{tab:experimentresults_Network}
\begin{tabular}{@{\hspace{5pt}}c@{\hspace{5pt}} @{\hspace{5pt}}c@{\hspace{5pt}} @{\hspace{5pt}}c@{\hspace{5pt}}
                @{\hspace{5pt}}c@{\hspace{5pt}} @{\hspace{5pt}}c@{\hspace{5pt}} @{\hspace{5pt}}c@{\hspace{5pt}}
                @{\hspace{5pt}}c@{\hspace{5pt}} @{\hspace{5pt}}c@{\hspace{5pt}}}
\hline
\multicolumn{2}{@{\hspace{5pt}}c@{\hspace{5pt}}}{} & \makecell[l]{Max-gap\\GP-UCB} & \makecell[l]{Max-gap\\sawtooth} & \makecell[l]{Random\\GP-UCB} & \makecell[l]{Random\\sawtooth} & \makecell[l]{Fixed-grid\\GP-UCB} & \makecell[l]{Fixed-grid\\sawtooth} \\[2pt]\hline
\multirow{5}{*}{h=10} & LB 
    & \textbf{151.180}\textsuperscript{(0.000)} 
    & \textbf{151.180} 
    & \textbf{151.180}\textsuperscript{(0.000)}
    & \textbf{151.180}\textsuperscript{(0.000)} 
    & 150.456 
    & 150.456 \\[2pt]
& Gap 
    & \textbf{0.010}\textsuperscript{(0.000)} 
    & \textbf{0.010} 
    & 1.348\textsuperscript{(0.232)} 
    & 1.846\textsuperscript{(0.357)} 
    & 104.631 
    & 141.734 \\[2pt]
& Time (s) 
    & 58.8\textsuperscript{(7.2)} 
    & 52.9
    & 3014.2\textsuperscript{(11.9)} 
    & 3012.2\textsuperscript{(7.0)} 
    & 0.7
    & 0.3\\[2pt]\cdashline{2-8}
& \makecell[c]{Worst Gap}  
    &  \textbf{0.010}
    &  \textbf{0.010}
    &  1.983
    &  2.365
    &  104.631
    &  141.734\\[2pt]\hline

\multirow{5}{*}{h=15} & LB 
    & \textbf{224.616}\textsuperscript{(0.000)} 
    & \textbf{224.616} 
    & 224.608\textsuperscript{(0.009)} 
    & 224.576\textsuperscript{(0.079)} 
    & 221.643 
    & 221.643 \\[2pt]
& Gap 
    & \textbf{0.025}\textsuperscript{(0.002)} 
    & 0.039 
    & 5.888\textsuperscript{(0.542)} 
    & 8.076\textsuperscript{(0.830)} 
    & 176.103 
    & 229.092\\[2pt]
& Time (s) 
    & 3011.7\textsuperscript{(13.5)} 
    & 3011.3
    & 3018.5\textsuperscript{(15.1)} 
    & 3012.9\textsuperscript{(10.4)} 
    & 1.2 
    & 0.5\\[2pt]\cdashline{2-8}
& \makecell[c]{Worst Gap}  
    &  \textbf{0.031}
    &  0.039
    &  6.535
    &  9.748
    &  176.103
    &  229.092\\[2pt]\hline

\multirow{5}{*}{h=20} & LB 
    & \textbf{298.149}\textsuperscript{(0.000)} 
    & \textbf{298.149} 
    & 298.101\textsuperscript{(0.032)} 
    & 298.053\textsuperscript{(0.084)} 
    & 293.685 
    & 293.685\\[2pt]
& Gap 
    & \textbf{0.164}\textsuperscript{(0.024)} 
    & 0.346
    & 12.579\textsuperscript{(1.042)} 
    & 15.545\textsuperscript{(0.995)} 
    & 248.946
    & 315.601\\[2pt]
& Time (s) 
    & 3033.7\textsuperscript{(20.6)} 
    & 3029.9
    & 3029.0\textsuperscript{(21.6)} 
    & 3011.9\textsuperscript{(11.1)} 
    & 1.7
    & 0.7\\[2pt]\cdashline{2-8}
& \makecell[c]{Worst Gap}  
    &  \textbf{0.216}
    &  0.346
    &  13.854
    &  17.795
    &  248.946
    &  315.601\\[2pt]\hline

\multirow{5}{*}{h=40} & LB 
    & \textbf{592.274}\textsuperscript{(0.003)} 
    & \textbf{592.274}
    & 591.878\textsuperscript{(0.150)} 
    & 591.719\textsuperscript{(0.229)} 
    & 582.704
    & 582.704\\[2pt]
& Gap 
    & \textbf{1.453}\textsuperscript{(0.018)} 
    & 2.419
    & 42.844\textsuperscript{(1.844)} 
    & 50.222\textsuperscript{(1.503)} 
    & 543.012 
    & 660.788\\[2pt]
& Time (s) 
    & 3030.6\textsuperscript{(36.2)} 
    & 3020.2
    & 3015.0\textsuperscript{(15.7)} 
    & 3020.2\textsuperscript{(13.1)} 
    & 3.2 
    & 1.5\\[2pt]\cdashline{2-8}
& \makecell[c]{Worst Gap}  
    &  \textbf{1.488}
    &  2.419
    &  45.825
    &  53.214
    &  543.012
    &  660.788\\[2pt]\hline

\end{tabular}
\end{table}

\begin{table}[htp]
\centering
\footnotesize
\renewcommand{\arraystretch}{1.0} % Adjust this to reduce or increase vertical spacing
\setlength{\tabcolsep}{4pt} % Adjust column separation if desired
\caption{Experiment results for the query problem. The target gap is $0.001$. Bolded numbers indicate the largest upper bounds and smallest gaps, values in parentheses are standard deviations. The fixed-grid method is unable to solve the problem within the limited 3000 seconds.}
\label{tab:experimentresults_Query}
\begin{tabular}{@{\hspace{5pt}}c@{\hspace{5pt}} @{\hspace{5pt}}c@{\hspace{5pt}} @{\hspace{5pt}}c@{\hspace{5pt}}
                @{\hspace{5pt}}c@{\hspace{5pt}} @{\hspace{5pt}}c@{\hspace{5pt}} @{\hspace{5pt}}c@{\hspace{5pt}}
                @{\hspace{5pt}}c@{\hspace{5pt}} @{\hspace{5pt}}c@{\hspace{5pt}}}
\hline
\multicolumn{2}{@{\hspace{5pt}}c@{\hspace{5pt}}}{} & \makecell[l]{Max-gap\\GP-UCB} & \makecell[l]{Max-gap\\sawtooth} & \makecell[l]{Random\\GP-UCB} & \makecell[l]{Random\\sawtooth} & \makecell[l]{Fixed-grid\\GP-UCB} & \makecell[l]{Fixed-grid\\sawtooth} \\[2pt]\hline
\multirow{5}{*}{h=10} & LB 
    & \textbf{30.021}\textsuperscript{(0.000)} 
    & \textbf{30.021} 
    & \textbf{30.021}\textsuperscript{(0.000)} 
    & \textbf{30.021}\textsuperscript{(0.000)} 
    & NA & NA \\[2pt]
& Gap 
    & \textbf{0.001}\textsuperscript{(0.000)} 
    & \textbf{0.001}
    & \textbf{0.001}\textsuperscript{(0.000)}
    & \textbf{0.001}\textsuperscript{(0.000)} 
    & NA & NA \\[2pt]
& Time (s) 
    & 5.7\textsuperscript{(0.6)} 
    & 80.7
    & 7.3\textsuperscript{(2.5)} 
    & 3030.8\textsuperscript{(7.9)} 
    & NA & NA \\[2pt]\cdashline{2-8}
& \makecell[c]{Worst Gap}  
    &  \textbf{0.001}
    &  \textbf{0.001}
    &  \textbf{0.001}
    &  \textbf{0.001}
    &  NA
    &  NA \\[2pt]\hline

\multirow{5}{*}{h=15} & LB 
    & \textbf{45.047}\textsuperscript{(0.000)} 
    & \textbf{45.047}
    & \textbf{45.047}\textsuperscript{(0.000)}
    & \textbf{45.047}\textsuperscript{(0.000)} 
    & NA & NA \\[2pt]
& Gap 
    & \textbf{0.001}\textsuperscript{(0.000)} 
    & 0.002
    & \textbf{0.001}\textsuperscript{(0.000)} 
    & 0.004\textsuperscript{(0.000)} 
    & NA & NA \\[2pt]
& Time (s) 
    & 8.6\textsuperscript{(0.6)} 
    & 3000.8
    & 30.6\textsuperscript{(8.4)}
    & 3036.2\textsuperscript{(13.9)} 
    & NA & NA \\[2pt]\cdashline{2-8}
& \makecell[c]{Worst Gap}  
    &  \textbf{0.001}
    &  0.002
    &  \textbf{0.001}
    &  0.004
    &  NA
    &  NA \\[2pt]\hline

\multirow{5}{*}{h=20} & LB 
    & \textbf{60.083}\textsuperscript{(0.000)} 
    & \textbf{60.083} 
    & \textbf{60.083}\textsuperscript{(0.000)}
    & \textbf{60.083}\textsuperscript{(0.000)} 
    & NA & NA \\[2pt]
& Gap 
    & \textbf{0.001}\textsuperscript{(0.000)} 
    & 0.004
    & \textbf{0.001}\textsuperscript{(0.000)} 
    & 0.007\textsuperscript{(0.000)} 
    & NA & NA \\[2pt]
& Time (s) 
    & 3007.6\textsuperscript{(4.8)} 
    & 3025.6
    & 206.2\textsuperscript{(72.9)}
    & 3039.5\textsuperscript{(21.0)} 
    & NA & NA \\[2pt]\cdashline{2-8}
& \makecell[c]{Worst Gap}  
    &  \textbf{0.001}
    &  0.004
    &  \textbf{0.001}
    &  0.007
    &  NA
    &  NA \\[2pt]\hline

\multirow{5}{*}{h=40} & LB 
    & \textbf{120.327}\textsuperscript{(0.000)} 
    & \textbf{120.327}
    & \textbf{120.327}\textsuperscript{(0.000)}
    & \textbf{120.327}\textsuperscript{(0.000)} 
    & NA & NA \\[2pt]
& Gap 
    & 0.017\textsuperscript{(0.001)} 
    & 0.026
    & \textbf{0.010}\textsuperscript{(0.001)}
    & 0.037\textsuperscript{(0.000)} 
    & NA & NA \\[2pt]
& Time (s) 
    & 3004.4\textsuperscript{(5.6)} 
    & 3049.5
    & 3018.8\textsuperscript{(32.9)}
    & 3027.4\textsuperscript{(30.6)} 
    & NA & NA \\[2pt]\cdashline{2-8}
& \makecell[c]{Worst Gap}  
    &  0.018
    &  0.026
    &  \textbf{0.010}
    &  0.037
    &  NA
    &  NA \\[2pt]\hline

\end{tabular}
\end{table}

\begin{table}[htp]
\centering
\footnotesize
\renewcommand{\arraystretch}{1.0} % Adjust this to reduce or increase vertical spacing
\setlength{\tabcolsep}{4pt} % adjust column separation if desired
\caption{Experiment results for the hallway problem. The target gap is $0.00001$. Bolded numbers indicate the largest upper bounds and smallest gaps, values in parentheses are standard deviations. The fixed-grid method is unable to solve the problem within the limited 3000 seconds.}
\label{tab:experimentresults_Hallway}
\begin{tabular}{@{\hspace{5pt}}c@{\hspace{5pt}} @{\hspace{5pt}}c@{\hspace{5pt}} @{\hspace{5pt}}c@{\hspace{5pt}}
                @{\hspace{5pt}}c@{\hspace{5pt}} @{\hspace{5pt}}c@{\hspace{5pt}} @{\hspace{5pt}}c@{\hspace{5pt}}
                @{\hspace{5pt}}c@{\hspace{5pt}} @{\hspace{5pt}}c@{\hspace{5pt}}}
\hline
\multicolumn{2}{@{\hspace{5pt}}c@{\hspace{5pt}}}{} & \makecell[l]{Max-gap\\GP-UCB} & \makecell[l]{Max-gap\\sawtooth} & \makecell[l]{Random\\GP-UCB} & \makecell[l]{Random\\sawtooth} & \makecell[l]{Fixed-grid\\GP-UCB} & \makecell[l]{Fixed-grid\\sawtooth} \\[2pt]\hline
\multirow{5}{*}{h=10} & LB 
    & 0.311\textsuperscript{(0.004)} 
    & 0.311
    & \textbf{0.327}\textsuperscript{(0.002)}
    & 0.326\textsuperscript{(0.003)} 
    & NA & NA \\[2pt]
& Gap 
    & \textbf{0.162}\textsuperscript{(0.000)} 
    & 0.168 
    & 0.188\textsuperscript{(0.005)}
    & 0.180\textsuperscript{(0.009)} 
    & NA & NA \\[2pt]
& Time (s) 
    & 3067.0\textsuperscript{(27.9)} 
    & 3034.5
    & 3054.3\textsuperscript{(36.8)}
    & 3065.6\textsuperscript{(47.5)} 
    & NA & NA \\[2pt]\cdashline{2-8}
& \makecell[c]{Worst Gap}  
    &  \textbf{0.162}
    &  0.168
    &  0.196
    &  0.196
    &  NA
    &  NA \\[2pt]\hline

\multirow{5}{*}{h=15} & LB 
    & \textbf{0.604}\textsuperscript{(0.001)} 
    & 0.587 
    & 0.591\textsuperscript{(0.012)} 
    & 0.587\textsuperscript{(0.011)} 
    & NA & NA \\[2pt]
& Gap 
    & \textbf{0.334}\textsuperscript{(0.001)} 
    & 0.360
    & 0.384\textsuperscript{(0.013)} 
    & 0.387\textsuperscript{(0.010)} 
    & NA & NA \\[2pt]
& Time (s) 
    & 2989.8\textsuperscript{(24.5)} 
    & 3082.8
    & 3086.2\textsuperscript{(57.0)}
    & 3062.6\textsuperscript{(33.9)} 
    & NA & NA \\[2pt]\cdashline{2-8}
& \makecell[c]{Worst Gap}  
    &  \textbf{0.336}
    &  0.360
    &  0.396
    &  0.391
    &  NA
    &  NA \\[2pt]\hline

\multirow{5}{*}{h=20} & LB 
    & \textbf{0.865}\textsuperscript{(0.001)} 
    & 0.858
    & 0.811\textsuperscript{(0.024)} 
    & 0.796\textsuperscript{(0.023)} 
    & NA & NA \\[2pt]
& Gap 
    & \textbf{0.507}\textsuperscript{(0.001)} 
    & 0.524
    & 0.604\textsuperscript{(0.030)} 
    & 0.615\textsuperscript{(0.026)} 
    & NA & NA \\[2pt]
& Time (s) 
    & 3022.0\textsuperscript{(25.1)} 
    & 2883.4
    & 3010.1\textsuperscript{(37.3)} 
    & 3063.4\textsuperscript{(41.6)} 
    & NA & NA \\[2pt]\cdashline{2-8}
& \makecell[c]{Worst Gap}  
    &  \textbf{0.508}
    &  0.524
    &  0.646
    &  0.635
    &  NA
    &  NA \\[2pt]\hline

\multirow{5}{*}{h=40} & LB 
    & \textbf{1.930}\textsuperscript{(0.000)} 
    & 1.912
    & 1.612\textsuperscript{(0.067)} 
    & 1.295\textsuperscript{(0.061)} 
    & NA & NA \\[2pt]
& Gap 
    & \textbf{1.220}\textsuperscript{(0.000)} 
    & 1.234
    & 1.451\textsuperscript{(0.070)} 
    & 1.510\textsuperscript{(0.061)} 
    & NA & NA \\[2pt]
& Time (s) 
    & 2871.8\textsuperscript{(18.9)} 
    & 3030.7
    & 3092.7\textsuperscript{(66.6)}  
    & 3069.2\textsuperscript{(21.0)} 
    & NA & NA \\[2pt]\cdashline{2-8}
& \makecell[c]{Worst Gap}  
    &  \textbf{1.220}
    &  1.234
    &  1.571
    &  1.638
    &  NA
    &  NA \\[2pt]\hline

\end{tabular}
\end{table}

\begin{table}[htp]
\centering
\footnotesize
\renewcommand{\arraystretch}{1.0} % Adjust this to reduce or increase vertical spacing
\setlength{\tabcolsep}{4pt} % Adjust column separation if desired
\caption{Experiment results for the aloha30 problem. The target gap is $0.01$. Bolded numbers indicate the largest upper bounds and smallest gaps, values in parentheses are standard deviations. The fixed-grid method is unable to solve the problem within the limited 3000 seconds.}
\label{tab:experimentresults_Aloha30}
\begin{tabular}{@{\hspace{5pt}}c@{\hspace{5pt}} @{\hspace{5pt}}l@{\hspace{5pt}} @{\hspace{5pt}}l@{\hspace{5pt}}
                @{\hspace{5pt}}l@{\hspace{5pt}} @{\hspace{5pt}}l@{\hspace{5pt}} @{\hspace{5pt}}l@{\hspace{5pt}}
                @{\hspace{5pt}}l@{\hspace{5pt}} @{\hspace{5pt}}l@{\hspace{5pt}}}
\hline
\multicolumn{2}{@{\hspace{5pt}}c@{\hspace{5pt}}}{} & \makecell[l]{Max-gap\\GP-UCB} & \makecell[l]{Max-gap\\sawtooth} & \makecell[l]{Random\\GP-UCB} & \makecell[l]{Random\\sawtooth} & \makecell[l]{Fixed-grid\\GP-UCB} & \makecell[l]{Fixed-grid\\sawtooth} \\[2pt]\hline
\multirow{5}{*}{h=10} & LB 
    & 122.527\textsuperscript{(0.000)} 
    & 122.501
    & \textbf{122.545}\textsuperscript{(0.006)} 
    & 122.517\textsuperscript{(0.016)} 
    & NA & NA \\[2pt]
& Gap 
    & \textbf{0.735}\textsuperscript{(0.000)} 
    & 0.736
    & 0.877\textsuperscript{(0.072)}
    & 0.963\textsuperscript{(0.121)} 
    & NA & NA \\[2pt]
& Time (s) 
    & 3042.3\textsuperscript{(6.2)} 
    & 2962.7
    & 3035.4\textsuperscript{(21.3)} 
    & 3004.7\textsuperscript{(22.7)} 
    & NA & NA \\[2pt]\cdashline{2-8}
& \makecell[c]{Worst Gap}  
    &  \textbf{0.735}
    &  0.736
    &  0.979
    &  1.177
    &  NA
    &  NA \\[2pt]\hline

\multirow{5}{*}{h=15} & LB 
    & 167.681\textsuperscript{(0.000)} 
    & 167.591
    & \textbf{167.784}\textsuperscript{(0.017)} 
    & 167.625\textsuperscript{(0.137)} 
    & NA & NA \\[2pt]
& Gap 
    & \textbf{1.836}\textsuperscript{(0.000)} 
    & 2.004
    & 2.026\textsuperscript{(0.132)} 
    & 2.460\textsuperscript{(0.229)} 
    & NA & NA \\[2pt]
& Time (s) 
    & 3050.9\textsuperscript{(21.2)} 
    & 2912.5 
    & 3034.1\textsuperscript{(29.2)}
    & 3101.9\textsuperscript{(77.6)} 
    & NA & NA \\[2pt]\cdashline{2-8}
& \makecell[c]{Worst Gap}  
    &  \textbf{1.836}
    &  2.004
    &  2.221
    &  3.017
    &  NA
    &  NA \\[2pt]\hline

\multirow{5}{*}{h=20} & LB 
    & 203.827\textsuperscript{(0.004)} 
    & 203.657
    & \textbf{204.090}\textsuperscript{(0.045)}
    & 203.765\textsuperscript{(0.203)} 
    & NA & NA \\[2pt]
& Gap 
    & \textbf{3.248}\textsuperscript{(0.007)} 
    & 4.003
    & 3.661\textsuperscript{(0.095)} 
    & 4.417\textsuperscript{(0.317)} 
    & NA & NA \\[2pt]
& Time (s) 
    & 3061.5\textsuperscript{(12.4)} 
    & 3070.2
    & 3073.9\textsuperscript{(32.9)} 
    & 3134.4\textsuperscript{(33.9)} 
    & NA & NA \\[2pt]\cdashline{2-8}
& \makecell[c]{Worst Gap}  
    &  \textbf{3.262}
    &  4.003
    &  3.804
    &  5.134
    &  NA
    &  NA \\[2pt]\hline

\multirow{5}{*}{h=40} & LB 
    & 282.129\textsuperscript{(0.000)} 
    & 282.039
    & \textbf{283.338}\textsuperscript{(0.249)}
    & 282.329\textsuperscript{(0.577)} 
    & NA & NA \\[2pt]
& Gap 
    & \textbf{10.428}\textsuperscript{(0.003)} 
    & 12.507
    & 11.164\textsuperscript{(0.632)}
    & 12.480\textsuperscript{(0.932)} 
    & NA & NA \\[2pt]
& Time (s) 
    & 2956.0\textsuperscript{(27.2)} 
    & 3086.8
    & 3084.9\textsuperscript{(45.8)}
    & 3163.7\textsuperscript{(39.4)} 
    & NA & NA \\[2pt]\cdashline{2-8}
& \makecell[c]{Worst Gap}  
    &  \textbf{10.435}
    &  12.507
    &  12.113
    &  14.570
    &  NA
    &  NA \\[2pt]\hline

\end{tabular}
\end{table}

\begin{figure}[htp]
  \centering
  \includegraphics[width=\linewidth]{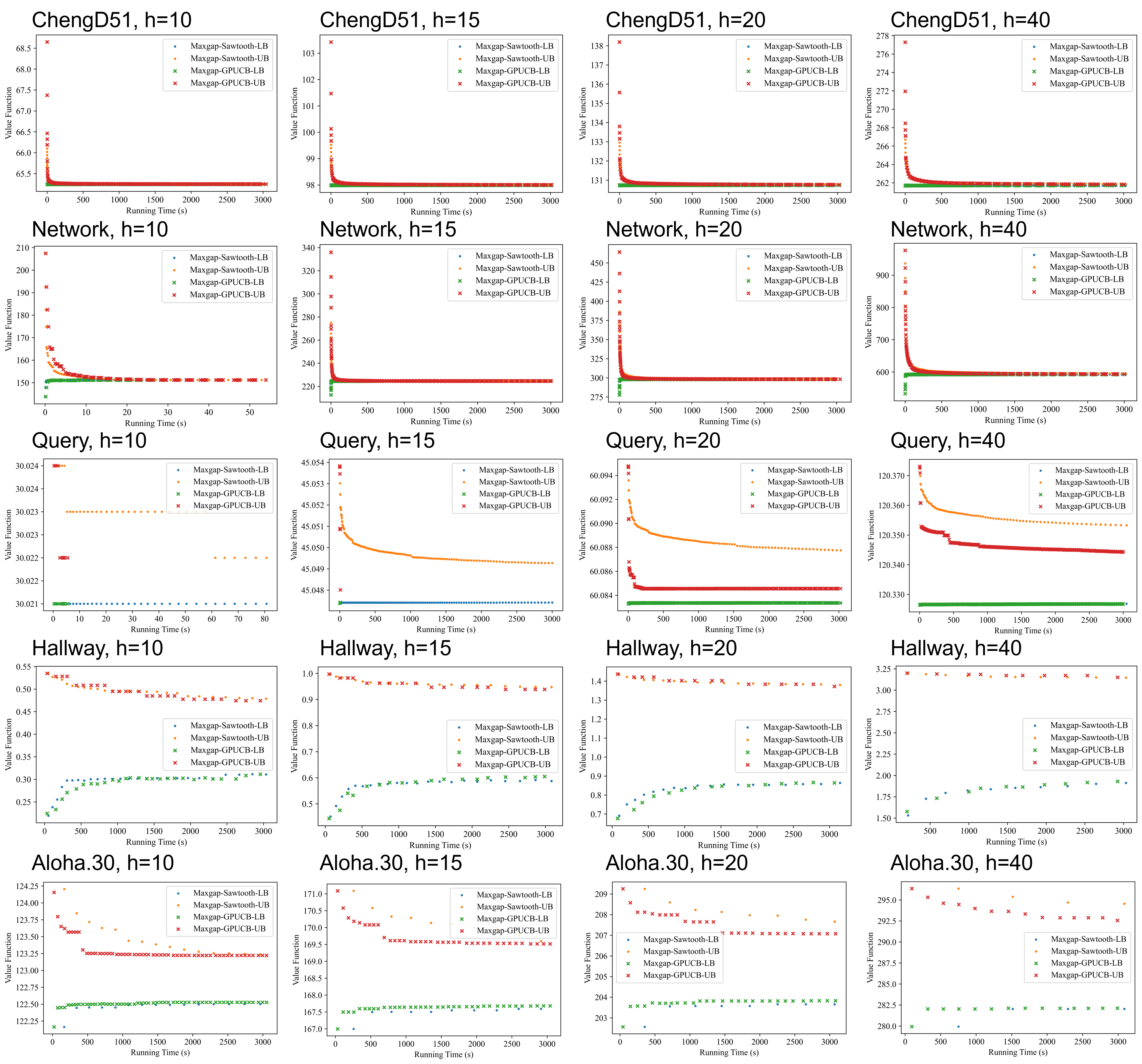}
  \caption{Upper and lower bounds of the value function for the five examples. Dots represent upper and lower bounds from FiVI, while crosses indicate those from AUCB, color-coded as shown in the legend.}
  \label{fig:ValueFunction}
\end{figure}

\end{document}